\documentclass[final]{article}

\usepackage{amsmath}
\usepackage{amsfonts}
\usepackage{color}
\usepackage{subfig}
\usepackage{graphicx}
\usepackage[margin=1in,headsep=.2in]{geometry}
\usepackage{hyperref}
\hypersetup{
    colorlinks=true,
    linkcolor=blue,
    citecolor=blue,
    filecolor=magenta,
    urlcolor=cyan,
}

\usepackage{algorithmic,algorithm}
\usepackage{enumitem}
\usepackage{amssymb}
\usepackage{comment}

\providecommand{\keywords}[1]
{
  \small
  \textbf{\textit{Keywords---}} #1
}
\usepackage{enumitem}
\def\zerobold{\boldsymbol{0}}
\def\onebold{\boldsymbol{1}}
\def\Abold{\boldsymbol{A}}

\def\abold{\boldsymbol{a}}

\def\rbold{\boldsymbol{r}}

\newlist{Assumption}{enumerate}{1}
\setlist[Assumption]{label=A\arabic*}

\newcommand{\bmat}[1]{\begin{bmatrix}#1\end{bmatrix}} 
\newcommand{\pmat}[1]{\begin{pmatrix}#1\end{pmatrix}} 
\definecolor{Blue}{rgb}{0,0,1}
\definecolor{Red}{rgb}{1,0,0}
\definecolor{Green}{rgb}{0,1,0}
\newcommand{\YC}[1]{{\textcolor{Red}{(YC: #1)}}}

\newcommand{\red}[1]{{\textcolor{Red}{#1}}}
\newcommand{\blue}[1]{{\textcolor{Blue}{#1}}}

\newcommand{\overbar}[1]{\mkern 1.5mu\overline{\mkern-1.5mu#1\mkern-1.5mu}\mkern 1.5mu}

\newcommand{\NN}{\mathbb{N}}
\newcommand{\RR}[1]{\ensuremath{\mathbb{R}^{ #1 }}}
\newcommand{\mapto}{\rightarrow}

\newcommand{\Span}[1]{\mathrm{span}\{#1\}}
\newcommand{\Range}[1]{\mathrm{range}(#1)}

\newcommand{\bandwidth}{b}
\newcommand{\fontDiscrete}{\mathcal}
\newcommand{\augmentedsubspace}{{\fontDiscrete G}}
\newcommand{\augmentedsubspaceArg}[1]{\augmentedsubspace_{#1}}
\newcommand{\subspace}{{\fontDiscrete A}}
\newcommand{\subspaceArg}[1]{\subspace_{#1}}
\newcommand{\snapshots}{{{\boldsymbol A}}}

\newcommand{\snapshotsArgII}[2]{\snapshots_{{#1}}^{{#2}}}
\newcommand{\systemmat}{\boldsymbol{K}}
\newcommand{\stiffness}{\systemmat}
\newcommand{\optvarSymbol}{x}
\newcommand{\optvar}{\boldsymbol{\optvarSymbol}}
\newcommand{\optvarMat}{\boldsymbol{X}}
\newcommand{\optvarArg}[1]{\optvarSymbol_{(#1)}}
\newcommand{\dummyvarISymbol}{q}

\newcommand{\dummyvarIVec}{\boldsymbol{\dummyvarISymbol}}
\newcommand{\dummyvarIISymbol}{\nu}
\newcommand{\dummyvarIIArg}[1]{\dummyvarIISymbol_{#1}}
\newcommand{\dummyvarIIVec}{\boldsymbol{\dummyvarIISymbol}}
\newcommand{\inverseoptvarSymbol}{z}

\newcommand{\inverseoptvarVec}{\boldsymbol{\inverseoptvarSymbol}}
\newcommand{\inverseoptvarMat}{\boldsymbol{Z}}
\newcommand{\onevec}{\boldsymbol{e}}
\newcommand{\pcgdiagonalmat}{\boldsymbol{\Gamma}}
\newcommand{\pcgdiagonalmatArg}[1]{\pcgdiagonalmat_{#1}}
\newcommand{\pcgdiagonalmatArgII}[2]{\pcgdiagonalmatArg{#1}^{(#2)}}
\newcommand{\krylovbasismat}{\boldsymbol{W}}
\newcommand{\krylovbasismatArg}[1]{\krylovbasismat_{#1}}
\newcommand{\krylovbasismatArgII}[2]{\krylovbasismatArg{#1}^{(#2)}}
\newcommand{\krylovsubspace}{{\mathcal K}}
\newcommand{\krylovsubspaceArg}[1]{\krylovsubspace_{#1}}
\newcommand{\krylovsubspaceArgII}[2]{\krylovsubspaceArg{#1}^{(#2)}}
\newcommand{\linsysvarSymbol}{u}
\newcommand{\linsysvar}{\boldsymbol{\linsysvarSymbol}}
\newcommand{\linsysvarInitial}{\linsysvar^{(0)}}
\newcommand{\linsysvarArg}[1]{\linsysvar_{#1}}
\newcommand{\linsysvarRef}{\linsysvar^{\text{ref}}}
\newcommand{\linsysvarRefArg}[1]{\linsysvarRef_{#1}}
\newcommand{\linsysvarArgII}[2]{\linsysvarArg{#1}^{(#2)}}
\newcommand{\linsysvarsol}{\linsysvar^\star}
\newcommand{\linsysvarsolArg}[1]{\linsysvarsol_{#1}}
\newcommand{\approximatelinsysvarSymbol}{\tilde{u}}
\newcommand{\approximatelinsysvar}{\boldsymbol{\approximatelinsysvarSymbol}}
\newcommand{\approximatelinsysvarArg}[1]{\approximatelinsysvar_{#1}}
\newcommand{\reducedlinsysvarSymbol}{\hat{u}}
\newcommand{\reducedlinsysvar}{\boldsymbol{\reducedlinsysvarSymbol}}
\newcommand{\reducedlinsysvarArg}[1]{\reducedlinsysvar_{#1}}
\newcommand{\reducedlinsysvarArgII}[2]{\reducedlinsysvarArg{#1}^{(#2)}}
\newcommand{\linsysmat}{\boldsymbol{A}}
\newcommand{\linsysmatArg}[1]{\linsysmat_{#1}}
\newcommand{\reducedlinsysmatArg}[1]{\hat{\linsysmat}_{#1}}
\newcommand{\rhs}{\boldsymbol{b}}

\newcommand{\rhsArg}[1]{\rhs_{#1}}
\newcommand{\card}[1]{|#1|}
\newcommand{\obj}{f}
\newcommand{\constraint}{g}
\newcommand{\constraintArg}[1]{\constraint_{#1}}
\newcommand{\constraintVec}{\boldsymbol{\constraint}}
\newcommand{\nonlinearequality}{h}

\newcommand{\nonlinearequalityVec}{\boldsymbol{\nonlinearequality}}
\newcommand{\pderes}{\rbold}
\newcommand{\pderesArg}[1]{\pderes_{#1}}
\newcommand{\reducedpderesArg}[1]{\hat{\pderes}_{#1}}
\newcommand{\stateSymbol}{w}
\newcommand{\state}{\boldsymbol{\stateSymbol}}
\newcommand{\stateArg}[1]{\state_{#1}}
\newcommand{\dispSymbol}{\stateSymbol}
\newcommand{\disp}{\boldsymbol{\dispSymbol}}

\newcommand{\decisionSymbol}{\mu}
\newcommand{\decision}{\boldsymbol{\decisionSymbol}}

\newcommand{\densitySymbol}{\rho}
\newcommand{\densityArg}[1]{\densitySymbol_{#1}}
\newcommand{\volumefractionSymbol}{\decisionSymbol}
\newcommand{\volumefractionArg}[1]{\volumefractionSymbol_{#1}}
\newcommand{\volumefraction}{\boldsymbol{\volumefractionSymbol}}
\newcommand{\lagrangemultiplierSenSymbol}{\lambda}
\newcommand{\lagrangemultiplierSen}{\boldsymbol{\lagrangemultiplierSenSymbol}}

\newcommand{\homotopyparam}{\omega}
\newcommand{\lagrangemultiplierSymbol}{\eta}
\newcommand{\lagrangemultiplier}{\boldsymbol{\lagrangemultiplierSymbol}}

\newcommand{\penaltySIMP}{s}
\newcommand{\penaltyStress}{q}
\newcommand{\stressRelaxed}{T}
\newcommand{\stressSolid}{T_0}
\newcommand{\youngsModulus}{E}
\newcommand{\youngsModulusSolid}{E_0}
\newcommand{\volumeSymbol}{v}

\newcommand{\volumeArg}[1]{\volumeSymbol_{#1}}

\newcommand{\nconstraint}{m}
\newcommand{\ndecision}{N_d}
\newcommand{\nelement}{\ndecision}
\newcommand{\nstate}{N_s}
\newcommand{\noptvar}{N_{\optvarSymbol}}
\newcommand{\designdomain}{\Omega}
\newcommand{\designdomainClosure}{\overbar{\designdomain}}
\newcommand{\designdomainArg}[1]{\designdomain_{#1}}
\newcommand{\designdomainClosureArg}[1]{\designdomainClosure_{#1}}
\newcommand{\up}{u}

\newcommand{\natNo}{\NN} 
\newcommand{\nat}[1]{\natNo(#1)}

\newcommand{\basismatSymb}{\Phi}
\newcommand{\basismat}{\boldsymbol{\basismatSymb}}
\newcommand{\basismatArg}[1]{\basismat_{#1}}
\newcommand{\middlebasismat}{\overbar{\basismat}}

\newcommand{\solSymb}{x}
\newcommand{\sol}{\boldsymbol \solSymb}

\newcommand{\defeq}{:=}

\newcommand{\opcgbasisvecSymbol}{z}
\newcommand{\opcgbasisvec}{\boldsymbol{\opcgbasisvecSymbol}}
\newcommand{\pcgbasisvecSymbol}{p}
\newcommand{\pcgbasisvec}{\boldsymbol{\pcgbasisvecSymbol}}
\newcommand{\resref}{r_0}

\newcommand{\preconditionerSymbol}{P}
\newcommand{\preconditionermat}{\boldsymbol{\preconditionerSymbol}}
\newcommand{\preconditioner}{\preconditionermat^{-1}}
\newcommand{\preconditionerArg}[1]{\preconditioner_{#1}}
\newcommand{\tol}{\epsilon}

\newcommand{\abstol}{\tol_{\text{abs}}}
\newcommand{\svdtol}{\tol_\text{SVD}}
\newcommand{\qrtol}{\tol_\text{QR}}
\newcommand{\betanom}{\eta}
\newcommand{\nom}{\zeta}
\newcommand{\den}{\gamma}
\newcommand{\stepsize}{\alpha}
\newcommand{\pcgbasisvecscale}{\beta}
\newcommand{\maxit}{N_\text{maxit}}

\newcommand{\inputvec}{\boldsymbol c}

\newcommand{\rightSV}{\boldsymbol{\Psi}}
\newcommand{\middlerightSV}{\overbar{\rightSV}}
\newcommand{\singularvaluesSymbol}{s}
\newcommand{\singularvalues}{\boldsymbol \singularvaluesSymbol}
\newcommand{\singularvaluesArg}[1]{\singularvaluesSymbol_{#1}}

\newcommand{\projected}{\boldsymbol \ell}
\newcommand{\projerr}{p}

\newcommand{\diag}[1]{\text{diag}(#1)}
\newcommand{\assign}{\leftarrow}
\newcommand{\Qmat}{\overbar{\boldsymbol{Q}}}
\newcommand{\Rmat}{\boldsymbol R}
\newcommand{\qr}[1]{ QR\{#1\} }

\newcommand{\bad}[1]{\red{\bf#1}}
\newcommand{\good}[1]{\blue{\bf#1}}

\newcommand{\leftsingularmat}{\boldsymbol{U}}
\newcommand{\rightsingularmat}{\boldsymbol{V}}
\newcommand{\singularvalmat}{\boldsymbol{\Sigma}}
\newcommand{\middlesingularvalmat}{\overbar{\singularvalmat}}
\newcommand{\basisrank}{r}
\newcommand{\maxrank}{r_{\text{max}}}

\newcommand{\normKKT}{r_{\text{kkt}}}
\newcommand{\normROM}{r_{\text{rom}}}
\newcommand{\identity}[1]{\boldsymbol{I}_{#1}}
\newcommand{\thresholdROM}{\epsilon_{\text{rom}}}
\newcommand{\thresholdPCG}{\epsilon_{\text{pcg}}}
\newcommand{\thresholdPCGupper}{\thresholdPCG^{\text{upper}}}
\newcommand{\thresholdPCGlower}{\thresholdPCG^{\text{lower}}}
\newcommand{\scaleCUT}{\kappa_\text{cut}}
\newcommand{\scaleROM}{\kappa_\text{rom}}
\newcommand{\scalePCG}{\kappa_\text{pcg}}

\newcommand{\absthresholdPCG}{\varepsilon_\text{pcg}}
\newcommand{\updatedBasisVector}{\boldsymbol{j}}
\newcommand{\reducedinputvec}{\boldsymbol{\ell}}
\newcommand{\middlemat}{\boldsymbol{Q}}
\newcommand{\helmholtzradius}{r}
\newcommand{\ipoptStationarityScaling}{s_d}
\newcommand{\ipoptComplimentarityScaling}{s_c}
\newcommand{\ipoptScaleMax}{s_{\max}}

\title{Accelerating design optimization using reduced order models}

\author{Youngsoo Choi\footnote{Correspondence to: Lawrence Livermore National
Laboratory, 7000 East Ave, Livermore, CA 94550, USA. E-mail: choi15@llnl.gov}, Geoffrey Oxberry, 
Daniel White, Trenton Kirchdoerfer\vspace{6pt}
\\Lawrence Livermore National Laboratory
\footnote{Lawrence Livermore National Laboratory is operated by Lawrence Livermore
National Security, LLC, for the U.S. Department of Energy, National Nuclear
Security Administration under Contract DE-AC52-07NA27344 and
LLNL-JRNL-791183.}
}
\date{}

\begin{document}

\maketitle

\begin{abstract}
Although design optimization has shown its great power of automatizing the whole
  design process and providing an optimal design, using sophisticated
  computational models, its process can be formidable due to a computationally
  expensive large-scale linear system of equations to solve, associated with
  underlying physics models.  We introduce a general reduced order model-based
  design optimization acceleration approach that is applicable not only to
  design optimization problems, but also to any PDE-constrained optimization
  problems. The acceleration is achieved by two techniques: i) {\it allowing an
  inexact linear solve} and ii) {\it reducing the number of iterations in Krylov
  subspace iterative methods}. The choice between two techniques are made, based
  on how close a current design point to an optimal point.  The advantage of the
  acceleration approach is demonstrated in topology optimization examples,
  including both compliance minimization and stress-constrained problems, where
  it achieves a tremendous reduction and speed-up when a traditional
  preconditioner fails to achieve a considerable reduction in the number of
  linear solve iterations. 
\end{abstract}

\keywords{design optimization, reduced order models, incremental SVD, KKT conditions,
  Krylov subspace iterative method, topology optimization, compliance
  minimization, stress-constrained problem}

\section{Introduction}

  Design optimization is a powerful tool that enables an automatic process of
  obtaining an optimal design with the help of sophisticated computational
  models. It is widely used in industries, academia, laboratories for various
  applications, mainly including aerospace, structural, mechanical, and
  biomedical engineering.  Many researchers also try to expand its physics
  domain to fluids, acoustics, electromagnetic, and optics. Thanks to the recent
  advance in additive manufacturing, a complicated optimal design can be
  directly manufactured.  However, the design optimization process involves an
  expensive physics model solution process. The most expensive part 
  is a sequence of large-scale linear solves. Although a sparse
  linear system may arise, the size of the system hinders a rapid design
  process.  There have been many attempts to reduce the cost of large-scale
  linear solves and they can be grouped to two categories: i) {\it allowing
  inexact linear solves} and ii) {\it reducing the number of iterations in
  Krylov subspace methods}.
 
  To {\it allow inexact linear solves} in design optimization, Amir and
  his coauthors in \cite{amir2010efficient} studied how to utilize the inexact
  solution of linear systems in the optimization process by considering specific
  convergence criteria and applying sensitivity correction terms to compliance
  and compliant mechanism problems. However, in order to extend their approach
  to other applications, such as stress-constrained problems, appropriate
  convergence criteria and sensitivity correction terms need to be developed.
  Gogu in \cite{gogu2015improving} replaced high-fidelity linear solve with
  Reduced Order Models (ROMs) to accelerate the design optimization process.
  Basis was constructed on-the-fly by Gram-Schmidt process whenever the ROM
  residual was bigger than a tolerance set by a user. In order to set the
  tolerance relatively a large value, e.g., 0.1, the sensitivity correction
  terms were added to the ROM sensitivity as in \cite{amir2010efficient}.
  However, the sensitivity correction terms required linear solves with full
  order model size.  Those solves involved system matrices whose factorizations
  were already known if a direct solver was used. However, if linear systems are
  solved with iterative solvers, those factorizations are not available.
  Additionally, the ROM method was compared with an academic version of topology
  optimization algorithm written in MATLAB. Therefore, it is likely that the
  reported speed-up will be degraded when the ROM-based method is compared with
  the High-Performance Computing (HPC)-based Full Order Model (FOM) solves. 
  Yoon in \cite{yoon2010structural} used various model reduction techniques for
  frequency response problems. The reduction methods were the mode
  superposition, Ritz vector, and quasi-static Ritz vector. All those ROMs
  mentioned above do not consider the reduced basis from the proper orthogonal
  decomposition, which is known to provide an optimal basis, given a data. 

  To {\it reduce the number of iterations in Krylov subspace methods}, a
  preconditioner is necessary. Although this paper does not focus on
  preconditioner (instead we use an existing preconditioner), the importance of
  a good preconditioner should not be ignored.  An optimal and efficient
  preconditioner depends on each problem. Most widely used preconditioners
  include, but not limited to, Jacobi preconditioner, incomplete Cholesky
  factorization \cite{golub1996}, Schur complement-based ones
  \cite{choi2015practical,rees2010all}, and multi-grid methods
  \cite{barker2016fast}.  Another way of reducing the number of iterations in
  Krylov subspace methods is to use a recycling approach.  The recycling approach
  has been mainly developed in numerical linear algebra and optimization
  communities, but some of the approaches were applied in topology optimization
  problems. For example, Wang and his coauthors in \cite{wang2007large} used
  MINRES \cite{paige1975solution} with recycling to accelerate the solution
  process of both symmetric positive-definite and indefinite systems.  Scaling
  of stiffness was used to bring down the condition number of stiffness matrix.
  They used incomplete Cholesky factorization as a preconditioner.  However, the
  recycling subspace was taken from the solutions of previous linear solves, not
  from reduced basis of previous solution, resulting in a bigger recycling
  space.  Unlike Wang's work, Carlberg and his coauthors in
  \cite{carlberg2016krylov} and Nguyen and Chen in \cite{nguyen2018reduced} used
  reduced basis of previous solutions, which is a compact representation of
  previous solutions. Thus, they were able to keep the dimension of recycling
  subspace small and still reduce the number of iterations.

  We present a novel ROM-based design optimization algorithm. It is motivated by
  the following characteristics of the optimization process:  The gradient-based
  optimization solvers start with an initial design and explore the design space
  until it finds an optimal solution that satisfies the KKT conditions. In the
  beginning of the optimization process, the change of design variables is
  large, indicating that the linear solves at this stage do not have to be
  solved precisely. This is exactly the motivation for several gradient-based
  optimization algorithms that introduce various inexactness to accelerate
  the optimization process (e.g., see
  \cite{bellavia1998inexact,heinkenschloss2002analysis,byrd2008inexact}).  The
  inexactness allows us to replace FOM with a ROM in the beginning of the
  optimization process where ROM can provide an approximate solution much faster
  than the corresponding FOM.  On the other hand, as the optimization process
  gets near the end, the change in design variables is small and a precise
  solution of each linear solve is required for a guarantee of the convergence
  to a local optimum.  For this stage of the optimization process, we will use
  Krylov subspace methods with ROM-recycling space that is able to solve for an
  accurate enough solution with a reduced number of iterations. The combination
  of {\it the inexact ROMs} in the beginning and {\it the Krylov subspace
  methods with ROM-recycling} near the convergence can accelerate the overall
  topology optimization process and maximize the usage of ROMs.

  Some additional features and contributions of our method are detailed in
  the following list: 
  \begin{itemize}
    \item Incremental Singular Value Decomposition (SVD) and Gram-Schmidt
      orthogonalization are used to determine
      the reduced basis on-the-fly. 
    \item Conjugate gradient with ROM-based recycling and Algebraic Multi-Grid
      (AMG) preconditioner are used to reduce the number of linear solve
      iterations. 
    \item Our ROM-based design optimization algorithm quickly finds a local
      optimal design that satisfies a necessary optimality condition, i.e., the
      Karush-Kuhn-Tucker (KKT) conditions.
    \item The norm of the KKT conditions are used to determine if the ROM can
      replace the FOM or not. 
    \item Our method is applied to a broad range of numerical examples of
      structural topology optimization: compliance minimization with single and
      multiple load cases and stress-constrained optimization with von Mises
      stress criterion. 
    \item Both structured and unstructured meshes are tested in numerical
      experiments. Our method is able to accelerate the problems with both
      structured and unstructured meshes, but more so for the problems
      with unstructured mesh.  
    \item Our proposed method is general enough so that it is applicable not
      only to topology optimizations, but also to general design optimization
      and any PDE-constrained optimization problems.
    \item Fully parallel version is implemented in C++ production code,
      developed at LLNL.
  \end{itemize}

  The rest of the paper is organized in the following way:
  Section~\ref{sec:topopt} presents a topology optimization formulation as an
  example of Partial Differential Equation (PDE)-constrained design optimization
  and its solution methodology. 
  Section~\ref{sec:IPM} presents the interior-point method as an example of
  optimization solvers that uses the norm of the KKT conditions as stopping
  criteria.  
  Section~\ref{sec:reducedbasis} explains how to construct an optimal reduced
  basis efficiently, using incremental factorization algorithms.
  Section~\ref{sec:ROMaccel} shows two ways of utilizing the reduced basis to
  accelerate the linear system solve process. In   
  Section~\ref{sec:ROM}, the first one, i.e., a projection-based reduced order
  model ROM is introduced.  
  In Section~\ref{sec:PCGrecycling}, the second one, i.e., a preconditioned conjugate gradient
  method with ROM-recycling is explained. 
  Various numerical examples of structural topology optimization problems are
  shown to demonstrate the advantages of our method in
  Section~\ref{sec:numeric}. Section~\ref{sec:conclusion} concludes with summary
  and future directions.

\subsection{Notations}\label{sec:notations}
    Scalars are denoted by lowercase letters, e.g., $a$.
    Vectors are denoted by boldface lowercase letters, e.g., $\abold$, and its
    $i$-th element by lowercase letter subscripted with an index, e.g., $a_i$.
    Matrices are denoted by boldface uppercase letters, e.g., $\Abold$, and its
    $(i,j)$-th element by lowercase letter subscripted with two indices, e.g.,
    $a_{ij}$. The cardinality of a subspace, $\subspace$, is denoted as
    $\card{\subspace}$. The real number space is denoted as $\RR{}$ and positive real number
    space as $\RR{}_+$. The range space of $\Abold$ is denoted as
    $\Range{\Abold}$. The 2-norm of a vector is defined as $\| \abold \|_2
    \defeq \sqrt{\sum_i a_i^2}$, 1-norm of a vector is defined as $\| \abold
    \|_1 \defeq \sum_i |a_i|$, and the infinity norm of a vector is defined as $\|
    \abold \|_{\infty} \defeq \max_i |a_i| $.

\section{Topology optimization}\label{sec:topopt}
\begin{figure}[t]
    \begin{center}
      \subfloat[classical design
      process]{\label{fig:topopt-flow-chart1}\includegraphics[scale=0.3]{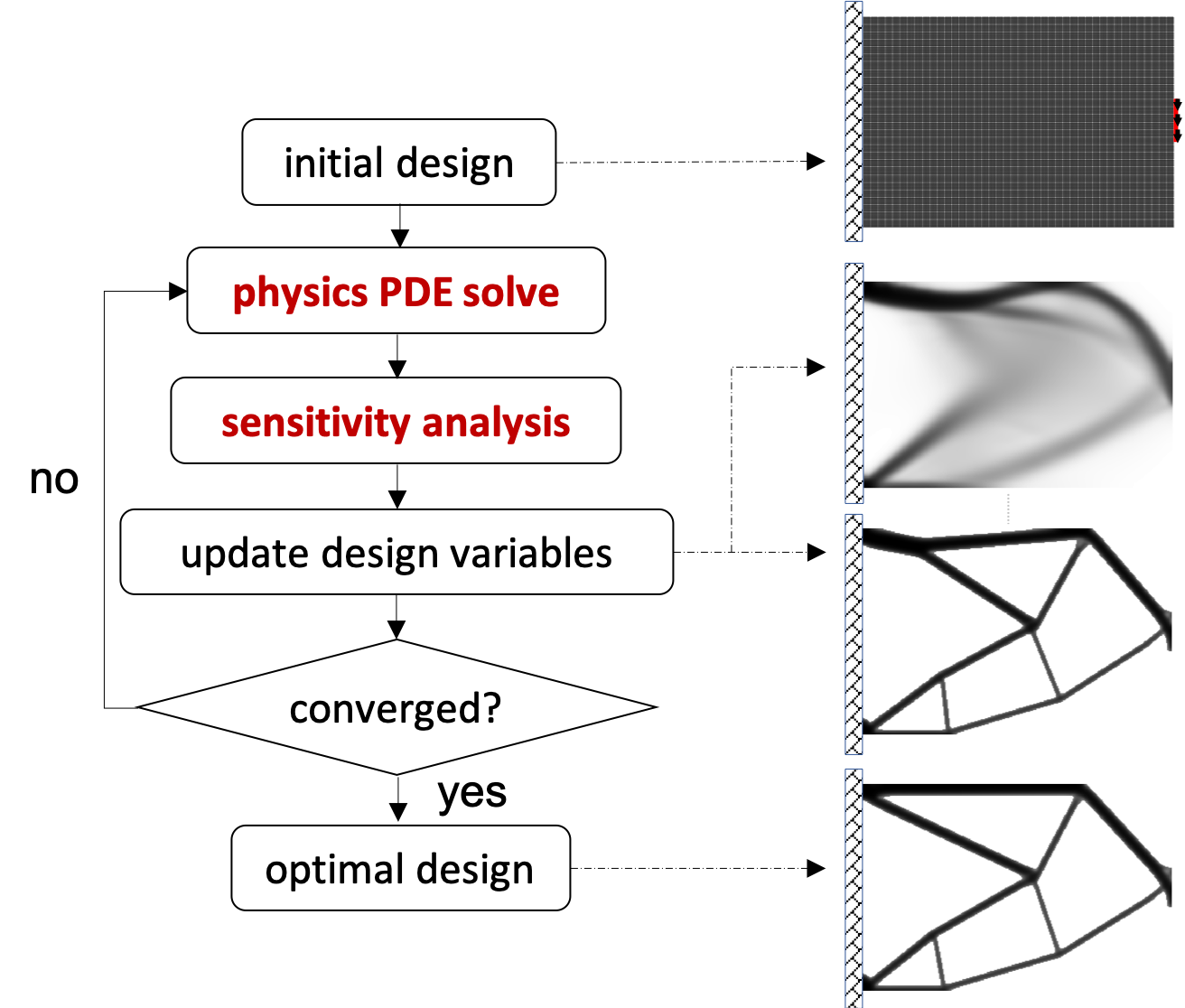}}
      \quad\quad\quad
      \subfloat[ROM-based design
      process]{\label{fig:topopt-flow-chart2}\includegraphics[scale=0.3]{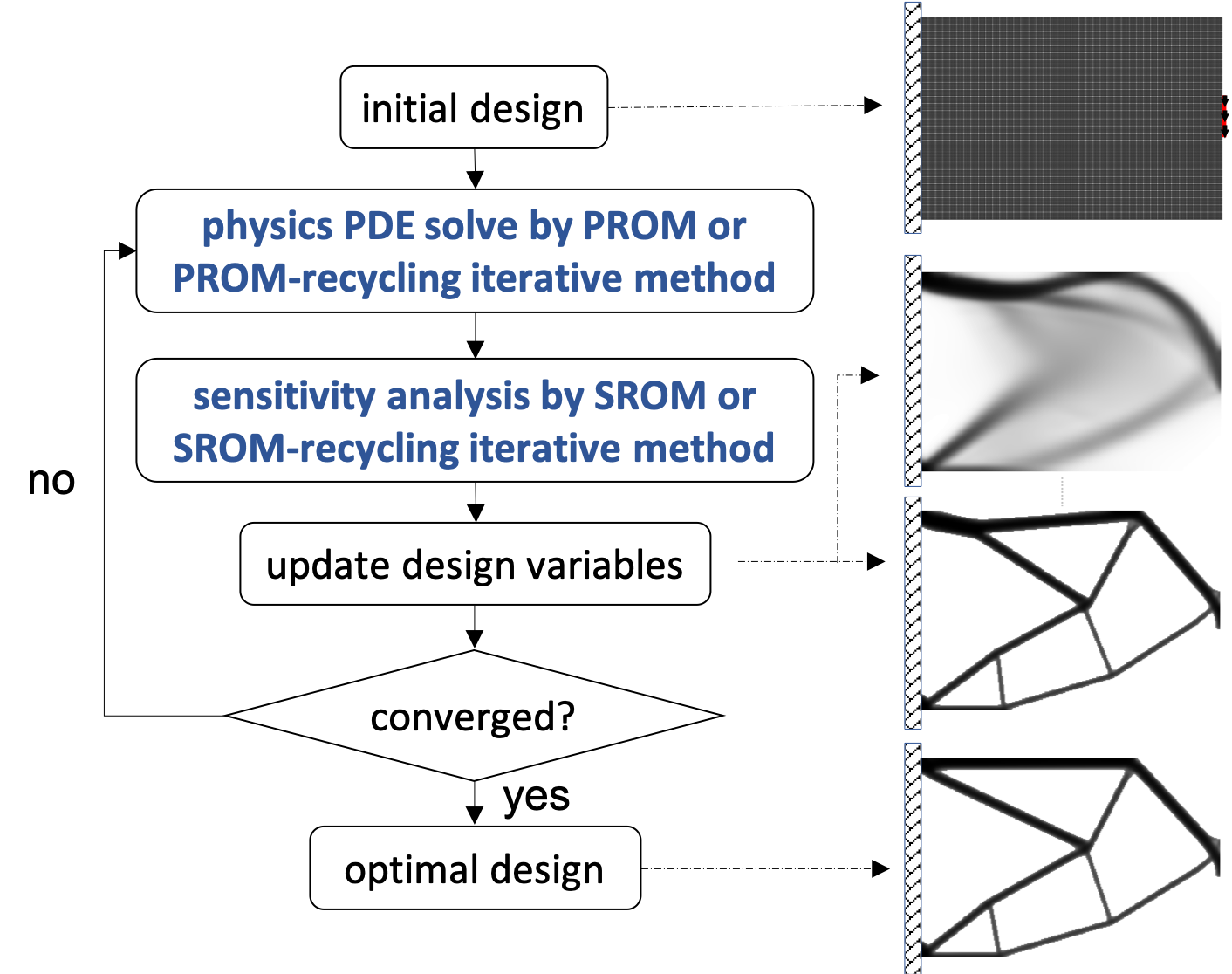}}
    \end{center}
    \caption{Comparison of classical and ROM-based design optimization flow
    charts in gradient-based optimization solvers. Basically, the ROM-based
    design process replaces the computationally expensive {\it PDE solves} and
    {\it sensitivity analysis} with either corresponding ROMs or ROM-recycling
    iterative methods. In (b), PROM refers to a ROM for {\it physics PDE
    solve} and SROM refers to a ROM for {\it sensitivity analysis}.}
    \label{fig:topopt-flow-chart}
  \end{figure}

  As an example of PDE-constrained design optimization, we consider a structural
  topology optimization problem. However, our method is general enough to be
  applicable to a broad class of design optimization problems.  A structural
  topology optimization finds the material distribution that minimizes an
  objective function, subject to $\nconstraint$ constraints. The material
  distribution is determined by the discretized volume fraction variables
  $\volumefraction\in\RR{\nelement}$ that can take any value between 0 (void)
  and 1 (solid material) for each non-overlapping element $\designdomainArg{i}$,
  $i\in\nat{\nelement}$, where $\nat{\nelement} := \{1,\ldots,\nelement\}$, in
  design domain $\designdomain$, and the number of elements $\nelement$ (i.e.,
  $\designdomainClosure = \cup_{i=1}^{\nelement} \designdomainClosureArg{i}$ and
  $\designdomainArg{i}\cap\designdomainArg{j} = []$ for $i\neq j$).
  This optimization problem can be mathematically formulated as:
  \begin{equation}\label{eq:topoptPb}
      \begin{aligned}
        & \underset{\volumefraction\in\RR{\nelement}}{\text{minimize}}  & & \obj(\disp(\volumefraction),\volumefraction)
        \\ & \text{subject to} & & \constraint_0(\volumefraction) = \sum_i
        \volumefractionArg{i}\volumeArg{i} - \volumeArg{\up} \leq 0
        \\ & & & \constraint_i (\disp(\volumefraction),\volumefraction) \leq 0,
        \quad i\in\nat{\nconstraint}
        \\ & & & \pderes(\state;\decision) = \zerobold
        \\ & & & \zerobold \leq \volumefraction \leq \onebold
      \end{aligned}
  \end{equation}
  where $f:\RR{\nstate}\times\RR{\nelement}\mapto\RR{}$ denotes an objective
  function, $\volumeArg{i}$ denotes the volume of element $\designdomainArg{i}$,
  $\volumeArg{\up}$ denotes the upper bound for the total volume of the
  material.  Optimization \eqref{eq:topoptPb} includes a volume constraint
  $\constraint_0 \leq 0$ and possibly $\nconstraint$ other nonlinear constraints
  $\constraint_i \leq 0$, $i\in\nat{\nelement}$.  Finally,
  $\disp:\RR{\nelement}\mapto\RR{\nstate}$ denotes a discretized displacement
  state vector function that depends implicitly on the volume fraction variables
  through the state equation residual function,
  $\pderes:\RR{\nstate}\times\RR{\ndecision}\mapto\RR{\nstate}$.  For example,
  the discretized PDE residual function for linear elasticity is defined as 
  \begin{equation}\label{eq:linearelasticity}
    \zerobold = \pderes(\state;\decision) := \systemmat(\decision)\state - \rhs(\decision),
  \end{equation}
  where $\stiffness:\RR{\nelement}\mapto\RR{\nstate\times\nstate}$ denotes a
  parameter dependent stiffness matrix (note $\stiffness =
  \frac{\partial\pderes}{\partial\state}$) and $\rhs:
  \RR{\nelement}\mapto\RR{\nstate}$ is a parameter-dependent right-hand-side
  vector for linear elasticity equation via finite element discretization (see
  Section 2 of \cite{hughes2012finite}). Note that the system of equations
  \eqref{eq:linearelasticity} need to be solved to evaluate $\obj$ and
  $\constraint_i$.  The solution process of \eqref{eq:linearelasticity} is
  labeled as the {\it physics PDE solve} in the topology optimization flow
  chart (see Figure~\ref{fig:topopt-flow-chart1}).

  Various quantities of interest can be considered as the objective, $\obj$, and
  constraint, $\constraint_i$, functions in the structural topology
  optimization.  For example, they include the compliance, i.e.,  $\state^T\rhs
  = \state^T\systemmat\state$ in discretized form \cite{bendsoe1989optimal},
  moment of inertia \cite{kang2005structural}, and various yield stress
  criteria, e.g., von Mises stress \cite{le2010stress} and Drucker-Prager for
  concrete materials \cite{luo2012topology}.

  In early days of topology optimization development, the checkerboard problem,
  i.e., patches of alternating void and material elements, was a major problem
  \cite{diaz1995checkerboard,jog1996stability,sigmund1998numerical}.  This
  problem is related to ensuring the well-posed and mesh-independent solutions.
  It can be addressed by introducing the various density filters: the cone
  filters \cite{bourdin2001filters,bruns2001topology} and the PDE filters
  \cite{lazarov2011filters,kawamoto2011heaviside}.  These filters act as
  smoothing and mixing tools of the irregular density variables, making
  individual density values dependent on neighboring elements.  Therefore, it
  should be noted that perimeter and slope constraints can also be used to avoid
  the checkerboard patterns.  Another possible filter is to apply the
  corresponding mass matrix to the density variables, which we call the mass
  filter. The minimum length scale of the mass matrix filter is imposed by the
  element size.  In this paper, we use either the mass matrix filter or the
  following Helmholtz type diffusion operator in \cite{lazarov2011filters}:
  \begin{equation}\label{eq:pde_filter}
    -\helmholtzradius^2 \Delta \densitySymbol + \densitySymbol = \volumefractionSymbol,
  \end{equation}
  where the filter radius, $r>0$ controls the minimum length scale and
  $\densitySymbol$ denotes a filtered volume fraction variable. Note that
  $\densitySymbol$ depends on $\volumefractionSymbol$ implicitly through
  \eqref{eq:pde_filter}.  Unfortunately, these filtering methods introduce grey
  areas between solid and void regions.  These grey areas can be minimized by
  reducing the filter radius. 

  The Solid Isotropic Material with Penalization (SIMP) or power-law approach is
  developed to ensure void-solid solutions by penalizing intermediate volume
  fraction variables \cite{zhou1991coc,mlejnek1992some,bendsoe1999material}. In
  the SIMP method, the volume fraction variable and the material property is
  related by
  \begin{equation}\label{eq:SIMP}
    \youngsModulus(\densityArg{i}) = \densityArg{i}^\penaltySIMP
    \youngsModulusSolid,
  \end{equation}
  where $\penaltySIMP\in(1,\infty)$ denotes the penalization parameter,
  $\youngsModulusSolid\in\RR{}$ denotes the Young's modulus of solid material, and
  finally, $\youngsModulus\in\RR{}$ denotes the penalized Young's modulus.
  A similar approach can be applied for stress-constrained topology optimization
  problems to address the {\it singularity} problem \cite{le2010stress}, i.e., 
  \begin{equation}\label{eq:relaxedStress}
    \stressRelaxed(\densityArg{i}) = \densityArg{i}^\penaltyStress
    \stressSolid,
  \end{equation}
  where $\penaltyStress\in(0,1)$ denotes the penalization parameter,
  $\stressSolid\in\RR{}$ denotes a stress quantity of solid material, and
  finally, $\stressRelaxed\in\RR{}$ denotes the relaxed stress. An alternative
  to the SIMP method is also available, e.g., the RAMP (Rational Approximation
  of Material Properties) method \cite{stolpe2001alternative} and the explicit
  penalization method \cite{allaire1993numerical,allaire1993topology}. As a 
  result, these penalized material properties are used to set up, for example,
  the state equations~\eqref{eq:linearelasticity} and various quantity of
  interests for $\obj$ and $\constraint_i$.

  A gradient-based optimization solver requires gradient
  $\frac{d\obj}{d\decision}\in\RR{\ndecision}$ and Jacobians
  $\frac{d\constraintArg{i}}{d\decision}\in\RR{\ndecision}$,
  $i\in\nat{\nconstraint}$. The {\it sensitivity analsys} in
  Figure~\ref{fig:topopt-flow-chart1}, is the step when those derivatives are
  computed.  They can be computed via chain rule. For example, the chain rule
  for the objective function gives:
  \begin{equation}\label{eq:chainrule}
    \frac{d\obj}{d\decision} = \frac{\partial\obj}{\partial\decision} +
    \frac{\partial\obj}{\partial\state} \frac{d\state}{d\decision},
  \end{equation}
  where $\frac{d\state}{d\decision}\in\RR{\nstate\times\ndecision}$ 
  can be obtained from the derivative of the linear PDE residual:
  \begin{align}
    \zerobold = \frac{d\pderes}{d\decision} 
       &= \frac{\partial\pderes}{\partial\decision} +
          \frac{\partial\pderes}{\partial\state} \frac{d\state}{d\decision},
          \label{eq:directsolve} \\
       &= \left ( \frac{d\systemmat}{d\decision}\state-\frac{d\rhs}{d\decision}
       \right ) + \systemmat(\decision)\frac{d\state}{d\decision}.
       \label{eq:detailedDirectSystem}
  \end{align}
  {\it The direct method} solves Eq.~\eqref{eq:directsolve} ``directly" for
  $\frac{d\state}{d\decision}$. Note that the direct method requires
  $\ndecision$ linear system solves for $\ndecision$ different right hand sides
  regardless of the number of quantity of interests, i.e., $\nconstraint+1$.
  Then, it substitutes $\frac{d\state}{d\decision}$ to \eqref{eq:chainrule} and
  obtains the gradient, $\frac{d\obj}{d\decision}$.  On the other hand, {\it the
  adjoint method} solves for Lagrange multipliers, $\lagrangemultiplierSen$, by
  solving the following linear adjoint system,
  \begin{equation}\label{eq:adjointsolve}
    \left (\frac{\partial\pderes}{\partial\state} \right
    )^T\lagrangemultiplierSen = \frac{\partial\obj}{\partial\state},
  \end{equation}
  then computes the gradient  
  \begin{equation}\label{eq:adjointgradient}
    \frac{d\obj}{d\decision} = \frac{\partial\obj}{\partial\decision} -
    \lagrangemultiplierSen\frac{\partial\pderes}{\partial\decision}.
  \end{equation}
  Note that the adjoint method requires one linear adjoint system solve for each
  quantity of interest (i.e., $\nconstraint+1$ linear adjoint system solves)
  regardless of the number of decision variables.  Therefore, one needs to use
  the adjoint method if $\nconstraint+1 < \ndecision$, while the direct method
  is preferred otherwise.  In density-based structural topology optimization,
  $\ndecision$ is proportional to the number of elements of spatial
  discretization, so it is likely to be larger than $\nconstraint+1$. Therefore,
  the adjoint method is preferred to the direct method in topology optimization.
  For more detailed review on topology optimization, please see these survey
  papers \cite{deaton2014survey,sigmund2013topology,sigmund1998numerical}.
  
  The gradient-based optimization algorithms need to solve
  a sequence of linear systems: 
  \begin{equation}\label{eq:multiplesolves}
    \linsysmatArg{k}\linsysvarArg{k} = \rhsArg{k},
  \end{equation}
  where the unknown variables are $\linsysvarArg{k} = \stateArg{k}$ for the PDE
  solve in \eqref{eq:linearelasticity}, $\linsysvarArg{k} =
  \frac{\partial\pderes_k}{\partial\state_k}$  for the direct sensitivity solve
  in \eqref{eq:directsolve}, and $\linsysvarArg{k} = \lagrangemultiplierSen_k$
  for the adjoint system in \eqref{eq:adjointsolve}. The subscript $k$ in
  $\state_k$, $\pderes_k$, and $\lagrangemultiplierSen_k$ indicates $k$th PDE,
  direct, and adjoint linear solves, respectively. These linear solves are
  necessary when objective, constraints, and their sensitivities need to be
  evaluated. Although a sparse linear system arises, the size and number of the
  linear system solution process hinders a rapid design process. These solves
  are the most expensive part of the optimization process.  For example, the
  computational time for the total optimization process for the wind turbine
  blade design problem in Section~\ref{sec:bladedesign} is $2.1$ hours, while the
  time for the linear solves is $1.7$ hours, taking about $81 \%$ of the total
  cost.  Therefore, developing a method of accelerating the sequence of
  large-scale linear solve is essential.  We achieve this by the ROM-based
  design optimization process described in Figure~\ref{fig:topopt-flow-chart2}.
  The ROM-based design optimization process enables to reduce $1.7$ hours of
  linear system solving time to $0.48$ hours for the same blade design problem,
  bringing down the cost of linear solves in the optimization process
  significantly.  Our method replaces the computationally expensive linear
  system solves occurred in {\it physics PDE solve} and {\it sensitivity
  analysis} indicated in Figure~\ref{fig:topopt-flow-chart1} with either {\it
  the inexact ROM solve} or {\it the ROM-recycling iterative solve}, e.g., see
  Figure~\ref{fig:topopt-flow-chart2}.  The choice between the inexact ROM or
  ROM-recycling iterative solve is determined by the norm of the KKT
  conditions.\footnote{The KKT conditions are the first-order optimality
  conditions if constraint qualifications are satisfied.  Many optimization
  algorithms use these conditions to look for a local minimum.  A sufficient
  condition for the local minimum is to have non-negative eigenvalues at the
  point.} In the following section, the KKT conditions in the context of the
  interior point method are derived.

\section{Interior point method}\label{sec:IPM}
  There are several gradient-based optimization algorithms that can solve
  Problem~\eqref{eq:topoptPb}.  For structural design optimizations, the Method
  of Moving Asymptotes (MMA) \cite{svanberg1987method}, the Optimality Criterion
  (OC)
  \cite{rozvany1989structural,rozvany1991coc,zhou1991coc,andreassen2011efficient},
  and CONLIN \cite{fleury1989conlin} are popular gradient-based algorithms.
  However, these algorithms are not robust in terms of the KKT conditions
  \cite{rojas2015benchmarking}. Unfortunately, this implies that these
  algorithms are not suitable to obtain a local minimum.  Our interest is in the
  optimization algorithms that can achieve the KKT 
  conditions.
  Such algorithms
  include, but not limited to,  the Sequential Quadratic Programming (SQP)
  \cite{boggs1995sequential}, the Interior Point Methods (IPMs)
  \cite{forsgren1998primal,forsgren2002interior,wachter2006implementation}, the
  trust region methods based on IPMs \cite{byrd2000trust}, and the augmented
  Lagrangian methods \cite{conn2013lancelot}. 
  We focus on the IPMs and its KKT conditions because we use IPOPT
  \cite{wachter2006implementation} 
  in our numerical experiment section. However, it is
  worthwhile to note that our method can be applied to any gradient-based
  optimization algorithms that find a KKT point.

 Introducing the following dummy variables and equations:
\begin{align}\label{eq:dummyVar}
  \dummyvarISymbol_0 = - \constraint_0,& \quad \dummyvarISymbol_0 \geq 0 \\ 
  \dummyvarISymbol_i = - \constraint_i,& \quad \dummyvarISymbol_i \geq 0,\quad
  i\in\nat{\nconstraint} \\ 
  \dummyvarIIArg{j} = 1-\volumefractionSymbol_j,& \quad \dummyvarIIArg{j} \geq 0,\quad
  j\in\nat{\nelement}, 
\end{align}
 and setting $\optvar^T = \pmat{\volumefraction^T & \dummyvarIVec^T & \dummyvarIIVec^T}$,
 Problem~\eqref{eq:topoptPb} can be transformed to the following general
 formulation:
\begin{equation}\label{eq:generalPb}
    \begin{aligned}
      & \underset{\optvar\in\RR{\noptvar}}{\text{minimize}}  & & \obj(\optvar)
      \\ & \text{subject to} & & \nonlinearequalityVec(\optvar) = \zerobold
      \\ & & & \optvar \geq \zerobold,
    \end{aligned}
\end{equation}
 where $\noptvar = 2\nelement+\nconstraint+1$ and
 $\nonlinearequalityVec(\optvar) = \constraintVec(\optvar) + \dummyvarIVec$. To
 solve \eqref{eq:generalPb}, IPM solves a sequence of barrier problems with the
 homotopy parameter $\homotopyparam\in\RR{}_+$ decreasing to zero:
\begin{equation}\label{eq:barrierPb}
    \begin{aligned}
      & \underset{\optvar\in\RR{\noptvar}}{\text{minimize}}  & & \obj(\optvar) -
      \homotopyparam \sum_{i=1}^{\noptvar} \ln(\optvarArg{i})
      \\ & \text{subject to} & & \nonlinearequalityVec(\optvar) = \zerobold.
    \end{aligned}
\end{equation}
Defining dummy variables, $\inverseoptvarSymbol_i = \homotopyparam /
\optvarSymbol_i$, $i\in\nat{\noptvar}$, the KKT conditions for
\eqref{eq:barrierPb} are 
 \begin{align}
   \nabla \obj(\optvar) + \nabla\nonlinearequalityVec(\optvar)\lagrangemultiplier -
   \inverseoptvarVec &= \zerobold \label{eq:stationarity} \\
   \nonlinearequalityVec(\optvar) &= \zerobold \label{eq:primalfeasibility} \\ 
   \optvarMat\inverseoptvarMat\onevec -  
   \homotopyparam \onevec &= \zerobold, \label{eq:complimentary}
 \end{align}
 where $\lagrangemultiplier\in\RR{\nconstraint+1}$ denotes a Lagrange
 multiplier vector, $\onevec \in \RR{\noptvar}$ denotes a vector with all the
 element values being one, and $\optvarMat$,
 $\inverseoptvarMat\in\RR{\noptvar\times\noptvar}$ are diagonal matrices whose
 diagonals are $\optvar$ and $\inverseoptvarVec$, respectively.  The interior
 point method follows the homotopy procedure where the barrier problem
 \eqref{eq:barrierPb} is solved for the unknowns,
 $(\optvar,\inverseoptvarVec,\lagrangemultiplier)$, with a fixed $\homotopyparam
 > 0$, then decreases the value of $\homotopyparam$, and solves the barrier
 problem with the previous solution as an initial guess. We repeat this process
 as the value of $\homotopyparam$ decreases to a small positive value close to
 zero. This process is justified by the fact that Eqs.~\eqref{eq:stationarity},
 \eqref{eq:primalfeasibility}, and \eqref{eq:complimentary} are the KKT
 conditions for Problem~\eqref{eq:generalPb} if $\homotopyparam = 0$ and
 $\optvar,\inverseoptvarVec \geq \zerobold$.\footnote{The KKT conditions of
 \eqref{eq:generalPb} include the stationarity
 condition~\eqref{eq:stationarity}, the primal feasibility
 conditions~\eqref{eq:primalfeasibility} and $\optvar\geq\zerobold$, the
 complementarity condition~\eqref{eq:complimentary} with $\homotopyparam = 0$,
 and duality condition, $\inverseoptvarVec\geq\zerobold$.} Therefore, when the
 homotopy procedure is complete, the solution to the barrier problem is a good
 approximate solution to Problem~\eqref{eq:generalPb}.

 One can include all the KKT equality conditions to define their norms: 
 the stationarity, primal feasibility, and complementarity 
 conditions. The IPOPT defines the norm of the KKT conditions in the following
 way:
 \begin{equation}\label{eq:KKTnorm}
   \normKKT = \max \left ( \frac{\| \nabla \obj(\optvar) +
   \nabla\nonlinearequalityVec(\optvar)\lagrangemultiplier - \inverseoptvarVec
   \|_\infty}{\ipoptStationarityScaling}, \| \nonlinearequalityVec(\optvar) \|_\infty,
   \frac{\|\optvarMat\inverseoptvarMat\onevec\|_\infty}{\ipoptComplimentarityScaling}
   \right),
 \end{equation}
 where $\ipoptStationarityScaling$, $\ipoptComplimentarityScaling \geq 1$ are
 defined as
 \begin{equation}\label{eq:ipoptScalings}
   \ipoptStationarityScaling = \max \left (\ipoptScaleMax,\frac{\|
   \lagrangemultiplier \|_1 + \| \inverseoptvarVec
   \|_1}{\nconstraint+\noptvar} \right ) / \ipoptScaleMax, \quad
   \ipoptComplimentarityScaling = \max \left ( \ipoptScaleMax,\frac{\|
   \inverseoptvarVec \|_1}{\noptvar}    \right ) / \ipoptScaleMax.
 \end{equation}
 This makes sure that a component of the optimality error is scaled, when the
 average value of the multipliers is larger than a fixed number
 $\ipoptScaleMax \geq 1$. We use $\ipoptScaleMax = 100$ as in IPOPT.
 The IPM software uses this norm or its variants
 to determine the convergence of the optimization process. The value of the KKT
 norm indicates how close the current point is to an optimal point.\footnote{A
 large $\normKKT$ indicates that the point is far from an optimal point. A
 small $\normKKT$ indicates the point is close to an optimal.} In ROM-based
 design optimization algorithm, FOM is completely replaced with ROM when the
 current point is far from an optimal point because the system does not need to
 be solved precisely.  Therefore, we will use $\normKKT$ as a measure to
 determine if the ROM replaces the corresponding FOM or not. 

\section{Reduced basis}\label{sec:reducedbasis}
  The topology optimization solves a sequence of linear system of equations
  \eqref{eq:multiplesolves} that generates a sequence of solutions,
  $\linsysvarArg{k}$.\footnote{Most gradient-based optimization algorithms find
  a local minimum.} The sequence of solutions converges to an optimal solution,
  $\linsysvar_\star$ if the problem is feasible. Additionally, as the
  optimization process is close to the end, the sequence of solutions does not
  change much, i.e., $\|\linsysvarArg{k} - \linsysvarArg{k+1}\|_2 \leq \epsilon$
  for $k > K > 0$ with a sufficiently large number $K$ and a small number
  $\epsilon$. Therefore, finding the solution $\linsysvarArg{k}$ within the
  subspace spanned by the previous $\ell$ solutions, $\subspaceArg{k-1}^\ell
  \defeq \Range{\snapshotsArgII{k-1}{\ell}}$, $\snapshotsArgII{k-1}{\ell}\defeq
  \bmat{\linsysvarArg{k-\ell},\ldots,\linsysvarArg{k-1}}
  \in\RR{\nstate\times\ell}$, must give a good approximation to
  $\linsysvarArg{k}$.  There are several ways of obtaining a basis $\basismat_k$
  for $\subspaceArg{k-1}^{\ell}$, but not limited to,
  \begin{itemize}
    \item a simple collection of previous $\ell$ solutions, i.e., $\basismat_k =
      \snapshotsArgII{k-1}{\ell}$
    \item an orthogonalized basis, e.g., via QR decomposition
    \item a Proper Orthogonal Decomposition (POD) basis.
  \end{itemize}
  The first choice is the simplest, but as the optimization process converges to
  an optimal solution, it generates a sequence of solution vectors that are
  almost linearly dependent. That causes the ill-condition of the reduced
  system.  Therefore, orthogonalization process, such as QR decomposition, is
  necessary, which is the second choice above.  The third choice above is
  motivated by POD.  The basis from POD is an optimally compressed
  representation of $\subspaceArg{k-1}$ in a sense that it minimizes the
  difference between the original snapshot matrix and the projected one onto the
  subspace spanned by the basis, $\basismat_k$:
  \begin{equation}\label{eq:POD}
    \begin{aligned}
      & \underset{\basismat_k\in\RR{\nstate\times\basisrank_k},
      \basismat_k^T\basismat_k
      =\identity{\basisrank_k}}{\text{minimize}} & & \left \|\snapshotsArgII{k-1}{\ell} -
      \basismat_k\basismat_k^T{\snapshotsArgII{k-1}{\ell}} \right \|_F^2, 
    \end{aligned}
  \end{equation}
  where $\|\cdot\|_F$ denotes the Frobenius norm and $\basisrank_k$ denotes the
  rank of the basis.  The solution of POD can be obtained by setting
  $\basismat_k = \leftsingularmat(:,1:\basisrank_k)$ in MATLAB notation, where
  $\leftsingularmat$ is the left singular matrix of the following thin Singular
  Value Decomposition (SVD):
 \begin{equation}\label{eq:SVD} 
   \snapshotsArgII{k-1}{\ell} = \leftsingularmat\singularvalmat\rightsingularmat^T,
 \end{equation} 
 where $\leftsingularmat\in\RR{\nstate\times\ell}$ and
 $\rightsingularmat\in\RR{\ell\times\ell}$ are orthogonal matrices and
 $\singularvalmat\in\RR{\ell\times\ell}$ is a diagonal matrix with singular
 values on its diagonals. SVD can order its modes from a most dominant mode to
 a least dominant mode. Thus, the first SVD basis vector is more important
 than the last SVD basis vector, making it easy to truncate and use only
 dominant modes in reduced basis.  POD is related to  principal component
 analysis in statistics \cite{hotelling1933analysis} and Karhunen--Lo\`{e}ve
 expansion \cite{loeve1955} in stochastic analysis.  Since the objective
 function in \eqref{eq:POD} does not change even though $\basismat_k$ is
 post-multiplied by an arbitrary $\basisrank_k\times\basisrank_k$ orthogonal
 matrix, the POD procedure seeks the optimal $\basisrank_k$-–dimensional
 subspace that captures the snapshots in the least-squares sense.  For more
 details on POD, we refer to
 \cite{berkooz1993proper,hinze2005proper,kunisch2002galerkin}.

 We will choose either the second or third choice above to generate our basis.
 However, it is computationally expensive to perform either QR or SVD of the
 snapshot matrix from scratch every time it is updated.  For example, the
 computational cost of SVD for
 $\snapshotsArgII{k-1}{\ell}\in\RR{\nstate\times\ell}$ is $O(\nstate^2\ell)$.
 Considering a large $\nstate$, this cost is too much.  Therefore, we use
 incremental algorithms where an efficient update to the previous decomposition
 is done when a new snapshot vector is added.  For the QR decomposition, the
 Gram-Schmidt orthogonalization perfectly fits into the incremental framework,
 i.e., see Algorithm~\ref{al:incrementalQR}. The inputs for the incremental QR
 in Algorithm~\ref{al:incrementalQR} include snapshot vector $\inputvec$,
 threshold for linear dependency $\qrtol$, index $k$, previous basis matrix
 $\basismat_{k-1}$, and maximum allowable rank of the basis matrix $\maxrank$.
 The incremental QR starts with $k=0$ with an empty basis matrix
 $\basismat_{-1}=[]$. For $k=0$, it simply normalizes the snapshot vector and
 set $\basismat_0$. For $k>0$, it applies the incremental QR update in
 Algorithm~\ref{al:incrementalQRupdate} as long as the rank of the basis matrix
 does not exceeds the maximum allowable rank $\maxrank$. If it exceeds the
 maximum allowable rank, then Algorithm~\ref{al:incrementalQR} throws away the
 first basis vector in $\basismat_{k-1}$. Then it applies the incremental QR
 update in Algorithm~\ref{al:incrementalQRupdate}. In the incremental QR update
 in Algorithm~\ref{al:incrementalQRupdate}, if snapshot vector $\inputvec$ can
 be spanned by the basis vectors in $\basismat_{k-1}$, then we set
 $\basismat_{k} = \basismat_{k-1}$. Otherwise, we include the effect of
 $\inputvec$ to $\basismat_{k}$.

  \begin{algorithm}[t]
    \caption{Incremental QR, $\basismat_{-1} = []$}\label{al:incrementalQR}
    $\basismat_k$ = \textbf{incrementalQR}($\inputvec$, $\qrtol$, $k$,
    $\basismat_{k-1}$,
    $\maxrank$) \\
    \textbf{Input:} $\inputvec$, $\qrtol$, $k$, $\basismat_{k-1}$, $\maxrank$\\
    \textbf{Output:} $\basismat_k$
      \begin{algorithmic}[1]
        \IF{$k=0$}
          \STATE $\basismat_0 \assign \inputvec/\| \inputvec \|$ 
        \ELSIF{$0<\basisrank_{k-1}\leq\maxrank$}
          \STATE $\basismat_k \assign \text{incrementalQRupdate}(\inputvec,
          \qrtol, \basismat_{k-1})$, i.e., apply
          Algorithm~\ref{al:incrementalQRupdate}
        \ELSE
          \STATE $\basismat_{k-1} \assign \basismat_{k-1}(:,2:\maxrank)$ 
          \STATE $\basismat_k \assign \text{incrementalQRupdate}(\inputvec,
          \qrtol, \basismat_{k-1})$, i.e., apply
          Algorithm~\ref{al:incrementalQRupdate}
        \ENDIF
   \end{algorithmic}
  \end{algorithm}

  \begin{algorithm}[t]
    \caption{Incremental QR update}\label{al:incrementalQRupdate}
    $\basismat_k$ = \textbf{incrementalQRupdate}($\inputvec$, $\qrtol$, $\basismat_{k-1}$)\\
    \textbf{Input:} $\inputvec$, $\qrtol$, $\basismat_{k-1}$\\
    \textbf{Output:} $\basismat_k$
      \begin{algorithmic}[1]
        \STATE $\updatedBasisVector \assign \inputvec-\basismat_k\basismat_k^T\inputvec$
        \IF{$\|\updatedBasisVector\| > \qrtol$}
          \STATE $\basismat_k \assign \bmat{\basismat_{k-1} & \updatedBasisVector/\|\updatedBasisVector\|}$  
        \ELSE
          \STATE $\basismat_k \assign \basismat_{k-1}$
        \ENDIF
   \end{algorithmic}
  \end{algorithm}
 
 For SVD, we use the incremental SVD in
 Algorithms~\ref{al:initializingIncrementalSVD} and \ref{al:incrementalSVD}.
 They were initially developed in \cite{oxberry2017limited}.  The algorithm is
 initialized with the {\it initializing incremental SVD} in
 Algorithm~\ref{al:initializingIncrementalSVD}. Then the following factorization
 is available for a rank--one update of the existing SVD
 \cite{brand2002incremental}: 

  \begin{align}
    \bmat{ \boldsymbol{\basismat}_{k-1} \singularvalmat_{k-1}
    \boldsymbol{\rightSV}_{k-1}^T & \inputvec}
    &= \bmat{ \boldsymbol{\basismat}_{k-1} & \left
    (\boldsymbol{I}-\boldsymbol{\basismat}_{k-1}\boldsymbol{\basismat}_{k-1}^T\right)
    \inputvec/\projerr} 
    \bmat{\singularvalmat_{k-1} &
    \boldsymbol{\basismat}_{k-1}^T\inputvec \\ \zerobold & \projerr}
    \bmat{ \boldsymbol{\rightSV}_{k-1} & \zerobold \\ \zerobold & 1}^T
    \label{eq:factorization1}\\
    &= \bmat{ \boldsymbol{\basismat}_{k-1} & \updatedBasisVector }
    \bmat{\singularvalmat_{k-1} & \reducedinputvec \\ \zerobold
    & \projerr}
    \bmat{ \boldsymbol{\rightSV}_{k-1} & \zerobold \\ \zerobold & 1}^T,
    \label{eq:factorization2}
  \end{align}
  where $\reducedinputvec = \boldsymbol{\basismat}_{k-1}^T\inputvec$ denotes a
  reduced coordinate of $\inputvec$, $\projerr = \| \inputvec -
  \basismat_{k-1}\reducedinputvec \|$ denotes the norm of the difference between
  $\inputvec$ and the projected one, and $\updatedBasisVector = \left
  (\inputvec-\boldsymbol{\basismat}_{k-1}\reducedinputvec\right)/\projerr$
  denotes a new orthogonal vector due to the incoming vector, $\inputvec$. Let
  \begin{equation}\label{eq:middleMat}
    \middlemat = \bmat{\singularvalmat_{k-1} & \reducedinputvec \\ \zerobold &
    \projerr}.
  \end{equation}
  The matrix,
  $\middlemat\in{\RR{(\basisrank_{k-1}+1)\times(\basisrank_{k-1}+1)}}$, is
  almost diagonal except for $\reducedinputvec$ in the upper right block and
  also its size is not in $O(\nstate)$. Thus, the SVD of $\middlemat$ is
  computationally fast. The cost is $O((\basisrank_{k-1}+1)^3)$, which is a lot
  cheaper than $O(\nstate^2\ell)$.
  Let the SVD of $\middlemat$ be
  \begin{equation}\label{eq:svdMiddleMat}
    \middlemat = \middlebasismat_{k-1}\middlesingularvalmat_{k-1}
    \middlerightSV_{k-1}, 
  \end{equation}
  where $\middlebasismat_{k-1} \in \RR{(\basisrank_{k-1}+1)\times(\basisrank_{k-1}+1)}$
  denotes the left singular matrix, $\middlesingularvalmat_{k-1} \in
  \RR{(\basisrank_{k-1}+1)\times(\basisrank_{k-1}+1)}$ denotes the singular
  value matrix, and $\middlerightSV_{k-1} \in
  \RR{(\basisrank_{k-1}+1)\times(\basisrank_{k-1}+1)}$ denotes the right
  singular matrix of $\middlemat$. Replacing $\middlemat$ in
  Eq.~\eqref{eq:factorization2} with \eqref{eq:svdMiddleMat} gives
  \begin{align}
    \bmat{ \boldsymbol{\basismat}_{k-1} \singularvalmat_{k-1}
    \boldsymbol{\rightSV}_{k-1}^T & \inputvec} &= 
    \bmat{ \boldsymbol{\basismat}_{k-1} & \updatedBasisVector }
    \middlebasismat_{k-1}\middlesingularvalmat_{k-1}
    \middlerightSV_{k-1}
    \bmat{ \boldsymbol{\rightSV}_{k-1} & \zerobold \\ \zerobold & 1}^T 
    \label{eq:rankOneUpdateSVD1}\\
    &= \boldsymbol{\basismat}_{k} \singularvalmat_{k}
    \boldsymbol{\rightSV}_{k}^T \label{eq:rankOneUpdateSVD2},
  \end{align}
  where $\boldsymbol{\basismat}_{k} = \bmat{ \boldsymbol{\basismat}_{k-1} &
  \updatedBasisVector } \middlebasismat_{k-1}\RR{\nstate\times(\basisrank_k)}$
  denotes the updated left singular matrix, $\singularvalmat_{k} =
  \middlesingularvalmat_{k-1} \in \RR{\basisrank_k\times\basisrank_k}$ denotes
  the updated singular value matrix, and $\boldsymbol{\rightSV}_{k} = \bmat{
    \boldsymbol{\rightSV}_{k-1} & \zerobold \\ \zerobold &
    1}\middlerightSV_{k-1} \in \RR{\basisrank_k\times\basisrank_k}$ denotes the
  updated right singular matrix.
  
  Algorithm~\ref{al:incrementalSVD} checks if $\inputvec$ is numerically
  linearly dependent on the current basis vectors. If $\projerr < \svdtol$, then
  we consider it is linearly dependent. Thus, we set $\projerr=0$ in
  $\middlemat$, i.e., Line 9 of Algorithm~\ref{al:incrementalSVD}. Then we only
  update the first $\basisrank_{k-1}$ components of the singular matrices in
  Line 14 of Algorithm~\ref{al:incrementalSVD}.  Although the orthogonality of
  the updated basis matrix, $\basismat_{k}$, must be guaranteed theoretically by
  the product of two orthogonal matrices in Line 14 or 16 of
  Algorithm~\ref{al:incrementalSVD}, it is not guaranteed numerically. Thus, we
  heuristically check the orthogonality in Lines 18-21 of
  Algorithm~\ref{al:incrementalSVD} by checking the inner product of the first
  and last columns of $\basismat_k$. If the orthogonality fails, then we
  orthogonalize them by the QR factorization.  Also, we limit the dimension of
  the basis to be less than or equal to $\maxrank$ because it is not necessary
  to include data far away from the current point of the optimization process in
  the reduced basis. See Line 1 of Algorithm~\ref{al:incrementalSVD}.

  \begin{algorithm}[t]
    \caption{Initializing incremental SVD}\label{al:initializingIncrementalSVD}
    [$\basismat_k$, $\singularvalues$, $\rightSV$] =
    \textbf{initializingIncrementalSVD}($\inputvec$, $\svdtol$, $k$) \\
    \textbf{Input:} $\inputvec$, $\svdtol$, $k$\\
    \textbf{Output:} $\basismat_k$, $\singularvalues$, $\rightSV$
      \begin{algorithmic}[1]
        \IF{$\|\inputvec\| > \svdtol$}
          \STATE $\singularvalues \assign  \bmat{\|\inputvec\|}$,
          $\basismat_k \assign \inputvec/ \singularvaluesArg{1}$, and
          $\rightSV \assign \bmat{1}$
        \ELSE
          \STATE $\singularvalues \assign []$,
          $\basismat_k \assign []$, and 
          $\rightSV \assign []$
        \ENDIF
   \end{algorithmic}
  \end{algorithm}

  \begin{algorithm}[ht]
    \caption{Incremental SVD, $\basismat_{-1} = [ ]$}\label{al:incrementalSVD}
    [$\basismat_k$, $\singularvalues$, $\rightSV$] =
    \textbf{incrementalSVD}($\inputvec$, $\svdtol$, $\basismat_{k-1}$,
    $\singularvalues$, $\rightSV$, $k$)\\
    \textbf{Input:} $\inputvec$, $\svdtol$, $\basismat_{k-1}$, $\singularvalues$,
    $\rightSV$, $k$\\
    \textbf{Output:} $\basismat_k$, $\singularvalues$, $\rightSV$
      \begin{algorithmic}[1]
        \IF{$\basisrank_{k-1}=0$ or $\basisrank_{k-1}=\maxrank$}
          \STATE [$\basismat_k$, $\singularvalues$, $\rightSV$] =
          initializingIncrementalSVD($\inputvec$, $\svdtol$, $k$), i.e., apply
          Algorithm~\ref{al:initializingIncrementalSVD} 
        \ENDIF
        \STATE $\boldsymbol \projected \assign \basismat_{k-1}^T\inputvec$ 
        \STATE $\projerr \assign \sqrt{\inputvec^T\inputvec-\projected^T\projected}$
        \STATE $\updatedBasisVector \assign (\inputvec - \basismat_{k-1}\projected)/\projerr$
        \STATE $\middlemat \assign \bmat{\diag{s} & \projected \\ \zerobold & \projerr}$ 
        \IF{$\projerr < \svdtol$}
          \STATE $\middlemat_{\text{end},\text{end}} \assign 0$
        \ENDIF
        \STATE
        $[\middlebasismat_{k-1},\middlesingularvalmat_{k-1},\middlerightSV_{k-1}]
        \assign \text{SVD}(\middlemat)$
        \STATE
        \IF{$\projerr<\svdtol$}
          \STATE $\basismat_{k} \assign
          \basismat_{k-1}{\middlebasismat_{k-1}}_{1:\basisrank_{k-1},1:\basisrank_{k-1}}$,
          \quad
          $\singularvalues_k \assign
          \diag{\middlesingularvalmat_{1:\basisrank_{k-1},1:\basisrank_{k-1}}}$,\quad
          and $\rightSV_k \assign {\middlerightSV_{k-1}}_{1:\basisrank_{k-1},:}$
        \ELSE
          \STATE $\basismat_{k} \assign \bmat{\basismat_{k-1} &
          \updatedBasisVector}\middlebasismat_{k-1}$, \quad $\singularvalues_k
          \assign \diag{\middlesingularvalmat_{k-1}}$, \quad and $\rightSV_k
          \assign \bmat{\rightSV_k & \zerobold \\ \zerobold &
          1}\middlerightSV_{k-1}$
        \ENDIF
        \IF{${\basismat_{k}}^T_{:,1} {\basismat_k}_{:,\text{end}} >
        \min\{\svdtol,\tol\cdot m_{\Phi}\}$}
          \STATE $[\Qmat,\Rmat] \assign \qr{\basismat_k}$ 
          \STATE $\basismat_k \assign \Qmat$
        \ENDIF
   \end{algorithmic}
  \end{algorithm}

  Once the basis matrix $\basismat_k$ is constructed, it can be used to
  construct ROM or recycling subspace of Krylov iterative methods.
  Specifically, Section~\ref{sec:ROMaccel} illustrates the flow chart of the
  ROM-based linear system acceleration scheme and Section~\ref{sec:ROM} explains
  how $\basismat_k$ is used to construct a projection-based ROM.
  Section~\ref{sec:PCGrecycling} shows the preconditioned conjugate gradient
  method with ROM-recycling. 

 \section{ROM-based linear system acceleration scheme}\label{sec:ROMaccel}
\begin{figure}[t]
    \begin{center}
      \includegraphics[scale=0.4]{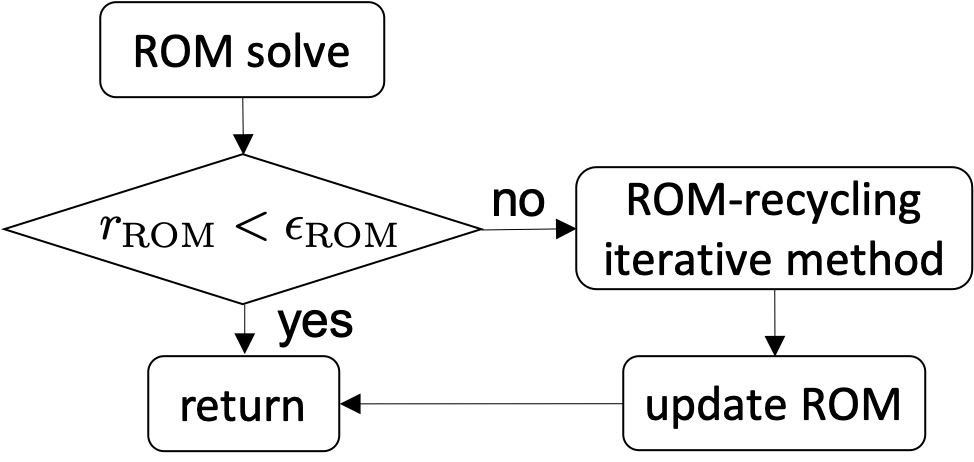}
    \end{center}
    \caption{Flow chart of the ROM-based linear system acceleration scheme}
    \label{fig:flowchart_linacceleration}
  \end{figure}

  To alleviate the cost of solving the system of linear
  equations~\eqref{eq:multiplesolves},\footnote{More precisely speaking, to
  alleviate the cost of solving the system of linear equations in
  \eqref{eq:linearelasticity} and either \eqref{eq:directsolve} or
  \eqref{eq:adjointsolve}} we apply the {\it ROM-based linear system
  acceleration scheme}, described in Figure~\ref{fig:flowchart_linacceleration}.
  The acceleration scheme starts with a ROM solve, which will be described in
  Section~\ref{sec:ROM}. If the ROM solution is good
  enough, then we use the ROM solution to update
  the design variables in Figure~\ref{fig:topopt-flow-chart2}.  Otherwise, it
  invokes the ROM-recycling iterative method, which gives a precise solution in
  a fast manner. The precise solution, in turn, is used as a new snapshot vector
  to update the current ROM within the incremental algorithms of
  Section~\ref{sec:reducedbasis} as well as the design variables in the
  optimization process of Figure~\ref{fig:topopt-flow-chart2}.
  Section~\ref{sec:ROM} shows how to build and solve ROMs. It also defines a
  residual norm of the ROM, $\normROM$, and a KKT conditions-related threshold,
  $\thresholdROM$. Section~\ref{sec:PCGrecycling} describes PCG with
  ROM-recycling method as an example of ROM-recycling iterative method.

  \subsection{Projection-based ROM}\label{sec:ROM} 

  Now we start to explain how the {\it ROM solve} is done in
  Figure~\ref{fig:flowchart_linacceleration}. The projection-based model
  reduction reduces the dimension of the system in \eqref{eq:multiplesolves} by
  reducing the number of the unknowns. For that purpose, the reduced basis
  $\basismat_k\in\RR{\nstate\times\basisrank_{k}}$ is used to approximate the
  solution variables as
 \begin{equation}\label{eq:approximatesol}
   \linsysvarArg{k} \approx \approximatelinsysvarArg{k}
   \defeq \linsysvarRefArg{k} + \basismat_k\reducedlinsysvarArg{k},
 \end{equation}
 where $\reducedlinsysvarArg{k}\in\RR{\basisrank_k}$ denotes the reduced
 coordinates and $\linsysvarRefArg{k}$ denotes a reference solution vector.
 Choices for $\linsysvarRefArg{k}$ include zero vector and  initial condition
 for time dependent problems.  Substituting \eqref{eq:approximatesol} into
 \eqref{eq:multiplesolves} and applying Galerkin projection lead to the
 following reduced system of equation:
 \begin{equation}\label{eq:ROM}
   \reducedlinsysmatArg{k}\reducedlinsysvarArg{k} =
   \reducedpderesArg{k},
   \basismat_k^T\pderesArg{k},
 \end{equation}
 where $\reducedlinsysmatArg{k} \defeq \basismat_k^T\linsysmatArg{k}\basismat_k
 \in\RR{\basisrank_k\times\basisrank_k}$ denotes a reduced linear system
 operator and $\reducedpderesArg{k} \defeq \basismat_k^T\left ( \rhsArg{j} -
 \linsysmatArg{j}\linsysvarRefArg{K} \right ) \in\RR{\basisrank_k}$ denotes a
 reduced right-hand side vector. The costs of constructing
 $\reducedlinsysmatArg{k}$ and $\reducedpderesArg{k}$ are at most
 $O(2(\bandwidth_k+1)\basisrank_k\nstate+\basisrank_k^3)$ and
 $O(\nstate\basisrank_k)$, respectively, assuming a banded structure in
 $\linsysmatArg{k}$ with its bandwidth, $\bandwidth_k$. The construction of
 $\reducedlinsysmatArg{k}$ can be large depending on how large $\basisrank_k$
 and $\bandwidth_k$ are. There is a hyper-reduction technique available to
 reduce the construction cost of reduced operators although we do not apply it
 to our numerical examples.  For example, see \cite{chaturantabut2010nonlinear,
 drmac2016new, choi2018sns} for various hyper-reduction techniques. 

 The cost for the solution process of \eqref{eq:ROM} is $O(\basisrank_k^3)$ for
 general methods.  The solution to Eq.~\eqref{eq:ROM} can be re-substituted to
 \eqref{eq:approximatesol} to recover the full size solution,
 $\approximatelinsysvarArg{k}$.  Since it is an approximation, the residual
 function value of $\approximatelinsysvarArg{k}$ is most likely to be non-zero,
 i.e., $\rhsArg{k} - \linsysmatArg{k}\approximatelinsysvarArg{k} \neq
 \zerobold$. We define the following relative residual norm of the ROM to
 measure the accuracy of $\approximatelinsysvarArg{k}$:
 \begin{equation}\label{eq:resnormROM}
   \normROM \defeq \frac{\| \rhsArg{k} - \linsysmatArg{k}\approximatelinsysvarArg{k}
   \|}{ \| \rhsArg{k} \|}. 
 \end{equation}
 Using this norm, we impose the condition, $\normROM<\thresholdROM$, for some
 threshold $\thresholdROM\in\RR{}_+$. The threshold is defined as 
 \begin{equation}\label{eq:ROMthreshold}
  \thresholdROM \defeq \scaleROM\cdot\normKKT,  
 \end{equation}
 where the reduction factor, $0<\scaleROM<1$, controls the tightness of the ROM
 accuracy and $\normKKT\in\RR{}_+$ denotes the norm of the KKT conditions.
 A smaller value of $\scaleROM$ expects a higher accuracy of the ROM.
 The norm of the KKT conditions, $\normKKT$, determines how far the current
 point is from the optimal design point. A large value of $\normKKT$ indicates
 that the current design variables are far from the optimal and a small value of
 $\normKKT$ indicates the other way around. Therefore, the ROM solutions are
 allowed if the current design variables are far from the optimal.  Otherwise,
 it is not allowed.

 Finally, if the condition in \eqref{eq:ROMthreshold} holds true, we consider
 that the ROM solution is good enough to be used as a solution of
 \eqref{eq:multiplesolves}. Otherwise, we invoke the ROM-recycling Krylov
 subspace method to compute a better solution than the ROM solution in a fast
 fashion. It is described in Section~\ref{sec:PCGrecycling}.

\subsection{ROM-recycling conjugate-gradient method}\label{sec:PCGrecycling}
  The system matrices arising from structural topology optimization problem,
  e.g., $\linsysmatArg{k}$ in Eq.~\eqref{eq:multiplesolves}, are usually sparse
  symmetric positive-definite matrices. Therefore, we consider the Conjugate
  Gradient (CG) method developed by Hestenes and Stiefel
  \cite{hestenes1952methods}.  However, we emphasize that our ROM-recycling
  approach can be applied to other iterative methods, such as MINRES
  \cite{paige1975solution, wang2007large} for symmetric indefinite systems and
  GCROT \cite{de1999truncation, parks2006recycling} and GMRES
  \cite{saad1997analysis} for more general nonsingular systems.  No matter what
  Krylov subspace iterative linear solver is used, a preconditioner needs to be
  applied to reduce the number of iterations for a large-scale problem.  Thus,
  we consider Preconditioned CG (PCG) and use Algebraic Multigrid Preconditioner
  (AMG) implemented in the scalable linear solver, HYPRE
  \cite{baker2012scaling}.

  The PCG method
  solves the following minimization problem at its $j$th iteration:
  \begin{equation}\label{eq:pcgIterateMinimization}
    \linsysvarArgII{k}{j} = \arg\min_{\linsysvar\in\linsysvarArgII{k}{0} +
    \krylovsubspaceArgII{k}{j}} \left \|
    \linsysvarsolArg{k} - \linsysvar \right \|_{\linsysmatArg{k}}
  \end{equation}
  where $\linsysvarArgII{k}{0}$ denotes a initial guess for the PCG process, 
  $\linsysvarsolArg{k}$ denotes the solution for the $k$th linear system,
  $\krylovsubspaceArgII{k}{j} \defeq \Span{\preconditionerArg{k}\pderesArg{k},
  \preconditionerArg{k}\linsysmatArg{k}\pderesArg{k}, \ldots,
  \preconditionerArg{k}\linsysmatArg{k}^{j-1}\pderesArg{k}}$ denotes the Krylov
  subspace at the $j$th PCG iteration of the $k$th linear system solve,
  $\pderesArg{k}\defeq \rhsArg{k} - \linsysmatArg{k}\linsysvarArgII{k}{0}$
  denotes the residual vector at the $k$th linear system solve,
  $\|\sol\|_{\linsysmatArg{k}}\defeq\sqrt{\sol\linsysmatArg{k}\sol}$ denotes the
  $\linsysmatArg{k}$-weighted norm, and $\preconditionerArg{k}$ denotes a
  preconditioner at $k$th linear system solve. The solution to
  Problem~\eqref{eq:pcgIterateMinimization} can be written as:
  \begin{equation}\label{eq:solutionPCG}
    \linsysvarArgII{k}{j} = \linsysvarArgII{k}{0} + \krylovbasismatArgII{k}{j}
    ({\krylovbasismatArgII{k}{j}}^T \linsysmatArg{k}
    \krylovbasismatArgII{k}{j})^{-1} {\krylovbasismatArgII{k}{j}}^T
    \pderesArg{k},
  \end{equation}
  where $\krylovbasismatArgII{k}{j}\in \RR{\nstate\times j}$ denotes a basis
  matrix for the Krylov subspace, $\krylovsubspaceArgII{k}{j}$. The basis
  matrix, $\krylovbasismatArgII{k}{j}$, satisfies
  $\linsysmatArg{k}$-orthogonality, i.e., 
    ${\krylovbasismatArgII{k}{j}}^T\linsysmatArg{k}\krylovbasismatArgII{k}{j} =
    \pcgdiagonalmatArgII{k}{j}$, where $\pcgdiagonalmatArgII{k}{j}$ is diagonal.
  Eq.~\eqref{eq:solutionPCG} can be viewed as the Galerkin projection,
  i.e., first solve for $\reducedlinsysvarArgII{k}{j}$ in the following reduced
  system:
  \begin{equation}\label{eq:pcgGalerkingProjection}
      {\krylovbasismatArgII{k}{j}}^T \linsysmatArg{k} \krylovbasismatArgII{k}{j}
      \reducedlinsysvarArgII{k}{j} =
      {\krylovbasismatArgII{k}{j}}^T\pderesArg{k}, 
  \end{equation}
  and then set $\linsysvarArgII{k}{j} = \linsysvarArgII{k}{0} +
  \krylovbasismatArgII{k}{j} \reducedlinsysvarArgII{k}{j}$. Note that this is
  identical procedure described in Eqs.~\eqref{eq:ROM} and
  \eqref{eq:approximatesol}.  
  
  The idea of recycling the Krylov subspace iterative linear solver comes from
  the context of solving a sequence of linear system of equations, e.g.,
  Eq.~\eqref{eq:multiplesolves}. It reuses data generated from the previous
  linear solves to reduce the iteration number for convergence of the current
  linear solve. That is, the recycling PCG algorithm introduces {\it augmented}
  subspace, $\augmentedsubspaceArg{k}$, in addition to the Krylov subspace,
  $\krylovbasismatArgII{k}{j}$, to jump start the PCG process, i.e., 
  \begin{equation}\label{eq:recyclingpcgIterateMinimization}
    \linsysvarArgII{k}{j} = \arg\min_{\linsysvar\in\linsysvarArgII{k}{0} +
    \augmentedsubspaceArg{k} + \krylovsubspaceArgII{k}{j}} \left \|
    \linsysvarsolArg{k} - \linsysvar \right \|_{\linsysmatArg{k}}.
  \end{equation}
  The augmented subspace $\augmentedsubspaceArg{k}$ is constructued, using the
  previous linear system solves \cite{wang2007large, o1980block,
  saad1987lanczos, erhel2000augmented, risler2000iterative}. It can be effective
  as long as the dimension of the subspace is not big. However, it is
  susceptible to the disadvantage of increasing the dimension of the augmented
  subspace as the number of linear system solves increases. Therefore, the
  truncation of the augmented subspace is needed. We accomplish the truncation,
  using the reduced basis, $\basismatArg{k}$, generated by the incremental
  algorithms in Section~\ref{sec:reducedbasis}. By compressing previous
  solutions in a reduced basis in an optimal way, i.e., in the POD sense, an
  important subspace information from previous solutions are kept within a
  small dimentional subspace. In summary, we define the augmented subspace,
  i.e., $\augmentedsubspaceArg{k} \defeq \Range{\basismatArg{k}}$. Then, the
  recycling PCG method generates the Krylov basis matrix,
  $\krylovbasismatArgII{k}{j}$, that is $\linsysmatArg{k}$-orthogonal to the
  augmented subspace,
  $\augmentedsubspaceArg{k}$, i.e.,
  \begin{equation}\label{eq:orthogonaltoAugmented}
    {\krylovbasismatArgII{k}{j}}^T\linsysmatArg{k}\basismatArg{k} =
    \zerobold,\quad\forall j 
  \end{equation}

  \begin{algorithm}[h]
    \caption{\red{ROM-recycling} PCG}\label{al:pcg}
    \YC{Redefine the notation here}\\
    $\linsysvar$ = \textbf{\red{ROMrecycling}PCG}($\linsysmat$, $\rhs$,
    $\linsysvarInitial$,
    $\thresholdPCG$, $\abstol$, $\maxit$, $\preconditioner$, \red{$\basismat$}) \\
    \textbf{Input:} $\linsysmat$, $\rhs$, $\linsysvar$, $\thresholdPCG$, $\abstol$,
    $\maxit$, $\preconditioner$, \red{$\basismat$} \\
    \textbf{Output:} $\linsysvar$
      \begin{algorithmic}[1]
        \STATE $\pderes \assign \rhs - \linsysmat\linsysvarInitial$ 
          \STATE $\pcgbasisvec \assign \preconditioner\pderes$
          \STATE \red{solve
          $\basismat^T\linsysmat\basismat\reducedlinsysvar =
          \basismat^T\linsysmat\pcgbasisvec$}
          \STATE \red{$\pcgbasisvec \assign \pcgbasisvec - \basismat\reducedlinsysvar$}
          \STATE $\resref \assign \max(\thresholdPCG^2\rhs^T\preconditioner\rhs,\abstol^2$)
          \STATE $\nom \assign \pderes^T\pcgbasisvec$
          \IF{$\nom \leq \resref$}
            \STATE converged
          \ENDIF
          \STATE $\opcgbasisvec \assign \linsysmat\pcgbasisvec$
          \STATE $\den \assign \opcgbasisvec^T\pcgbasisvec$
          \IF{$\den = 0$}
            \STATE fail to converge
          \ENDIF
          \FOR{ $j = 1,\ldots,\maxit$ }
            \STATE $\stepsize \assign \frac{\nom}{\den}$ 
            \STATE $\linsysvar \assign \linsysvar + \stepsize\pcgbasisvec$
            \STATE $\pderes \assign \pderes - \stepsize\linsysmat\pcgbasisvec$
            \STATE $\opcgbasisvec \assign \preconditioner\pderes$
            \STATE $\betanom \assign \pderes^T\opcgbasisvec$
            \IF{$\betanom < \resref$}
              \STATE converged
            \ENDIF
            \STATE $\pcgbasisvecscale \assign \frac{\betanom}{\nom}$
            \STATE $\pcgbasisvec \assign \opcgbasisvec + \pcgbasisvecscale\pcgbasisvec$

            \STATE \red{solve
            $\basismat^T\linsysmat\basismat\reducedlinsysvar =
            \basismat^T\linsysmat\opcgbasisvec$}
            \STATE \red{$\pcgbasisvec \assign \pcgbasisvec - \basismat\reducedlinsysvar$}

            \STATE $\opcgbasisvec \assign \linsysmat\pcgbasisvec$
            \STATE $\den \assign \pcgbasisvec^T\opcgbasisvec$
            \IF{ $\den \leq 0$ }
              \STATE not positive definite
            \ENDIF
            \STATE $\nom \assign \betanom$
          \ENDFOR
   \end{algorithmic}
  \end{algorithm}

  Algorithm \ref{al:pcg} describes the ROM-recycling PCG method. For the
  brevity, we skip both the subscripts and superscripts of each variables. If
  the red blocks are omitted, then the algorithm falls into a usual PCG method.
  Note that the Krylov subspace basis vector, $\pcgbasisvec$, is modified by the
  solution from Galerkin projection, i.e., Lines 3--4 and 26--27 of
  Algorithm~\ref{al:pcg}. The simple subtraction of the Galerkin part from
  $\pcgbasisvec$ is to ensure the $\linsysmatArg{k}$-orthogonality expressed in
  Eq.~\eqref{eq:orthogonaltoAugmented}.

  We set the initial guess for the ROM-recycling PCG to be the solution of the
  {\it ROM solve} in Figure~\ref{fig:flowchart_linacceleration}, i.e.,
  $\linsysvarInitial = \approximatelinsysvarArg{k}$ of
  Eq.~\eqref{eq:approximatesol}. The motivation for this choice can be explained
  by the fact that $\approximatelinsysvarArg{k}$ of
  Eq.~\eqref{eq:approximatesol} is the solution of the following minimization
  problem:
  \begin{equation}\label{eq:romMinimization}
    \approximatelinsysvarArg{k} = \arg\min_{\linsysvar\in\linsysvarRefArg{k} +
    \augmentedsubspaceArg{k} } \left \| \linsysvarsolArg{k} - \linsysvar \right
    \|_{\linsysmatArg{k}}.
  \end{equation}
  Therefore, this choice of the initial guess makes the optimal starting point
  in the sense of Eq.~\eqref{eq:romMinimization}. 

  As in the projection-based ROM of Section~\ref{sec:ROM}, we set the PCG
  convergence threshold, $\thresholdPCG$, relative to the value of the norm of
  the KKT conditions, $\normKKT$. More specifically, we provide the following
  two options:
  \begin{enumerate}
    \item The first one sets $\thresholdPCG$ in the following way:
      \begin{equation}\label{eq:PCGthreshold}
        \thresholdPCG \defeq \min(\max(\scalePCG\cdot\normKKT,
        \thresholdPCGlower), \thresholdPCGupper),
      \end{equation}
      where the reduction factor, $0<\scalePCG<1$, controls the tightness of the
      PCG accuracy. A smaller value of $\scalePCG$ expects a higher accuracy of
      the PCG solution. Also, we include the norm of the KKT conditions,
      $\normKKT$, in the definition of $\thresholdPCG$ to determine how far the
      current point is from the optimal design point. A large value of
      $\normKKT$ indicates that the current design variables are far from the
      optimal and a small value of $\normKKT$ indicates the other way around.
      Therefore, the less precise PCG solutions are allowed if the current
      design variables are far from the optimal. Otherwise, the more accurate
      PCG solution must be computed.  Additionally, we introduce the two
      safeguard thresholds, $\thresholdPCGlower$ and $\thresholdPCGupper$, such
      as $\thresholdPCGlower < \thresholdPCGupper$. The lower bound of the PCG
      convergence threshold is set by $\thresholdPCGlower$, while
      $\thresholdPCGupper$ serves as the upper bound. This is to avoid the case
      when the KKT norm is either too large or too low. 
      
    \item The second choice explicitly introduces a cut value, $\scaleCUT$, for
      the PCG convergence threshold. If $\normKKT > \scaleCUT$, then we set
      $\thresholdPCG = \max(\scalePCG\cdot\normKKT, \absthresholdPCG)$.
      Otherwise, we set $\thresholdPCG = \absthresholdPCG$. This implies that if
      the norm of the KKT conditions is greater than the cut value, then we
      allow the inexact solve by the PCG because the current point is far away
      from an optimal point and the linear solve does not need to be solved
      precisely. If $\normKKT \leq \scaleCUT$, we set the PCG convergence
      threshold to be a user-defined threshold, i.e., $\absthresholdPCG$.
      Therefore, the minimum threshold this option can take is
      $\absthresholdPCG$. 
  \end{enumerate}

\section{Numerical experiments}\label{sec:numeric}

  We provide numerical evidences of the advantages of our method in
  several numerical experiments with structural topology optimization.
  We consider two different types of topology optimization problems: 1)
  compliance minimization problem with mass constraint and 2) the mass
  minimization problems with stress constraint.

  For all the numerical examples considered herein, we use the second option of
  setting $\thresholdPCG$ described in Section~\ref{sec:PCGrecycling}. All the
  simulations use a number of processors in Quartz of Livermore Computing
  Center\footnote{https://hpc.llnl.gov/hardware/platforms/Quartz}. All the
  visualizations are made with
  VisIt\footnote{https://wci.llnl.gov/simulation/computer-codes/visit}.
  The default and our ROM-based topology optimization methods are compared.
  The default method follows the flow chart in
  Figure~\ref{fig:topopt-flow-chart1} where {\it physics PDE solve} and {\it
  sensitivity analysis} are solved by PCG without ROM-recycling, i.e., see
  Algorithm~\ref{al:pcg} without red blocks, but with AMG preconditioner from
  HYPRE. Zero initial guess is used in the PCG iterations and we set
  $\thresholdPCG = \absthresholdPCG$ for the default method. For the ROM-based
  topology optimization methods, we consider both incremental QR (see Algorithms
  \ref{al:incrementalQR} and \ref{al:incrementalQRupdate}) and SVD (see
  Algorithms \ref{al:initializingIncrementalSVD} and \ref{al:incrementalSVD}) to
  construct reduced bases. Both incremental algorithms are implemented in an
  open source library, libROM \cite{osti_1505575}.

  \subsection{Compliance minimization}\label{sec:compliance}
  Two 3D design problems are considered: i) cantilever beam, and ii) wind turbine
  blade. 

  \subsubsection{3D cantilever beam design}
\begin{figure}[th]
    \begin{center}
    \includegraphics[scale=0.30]{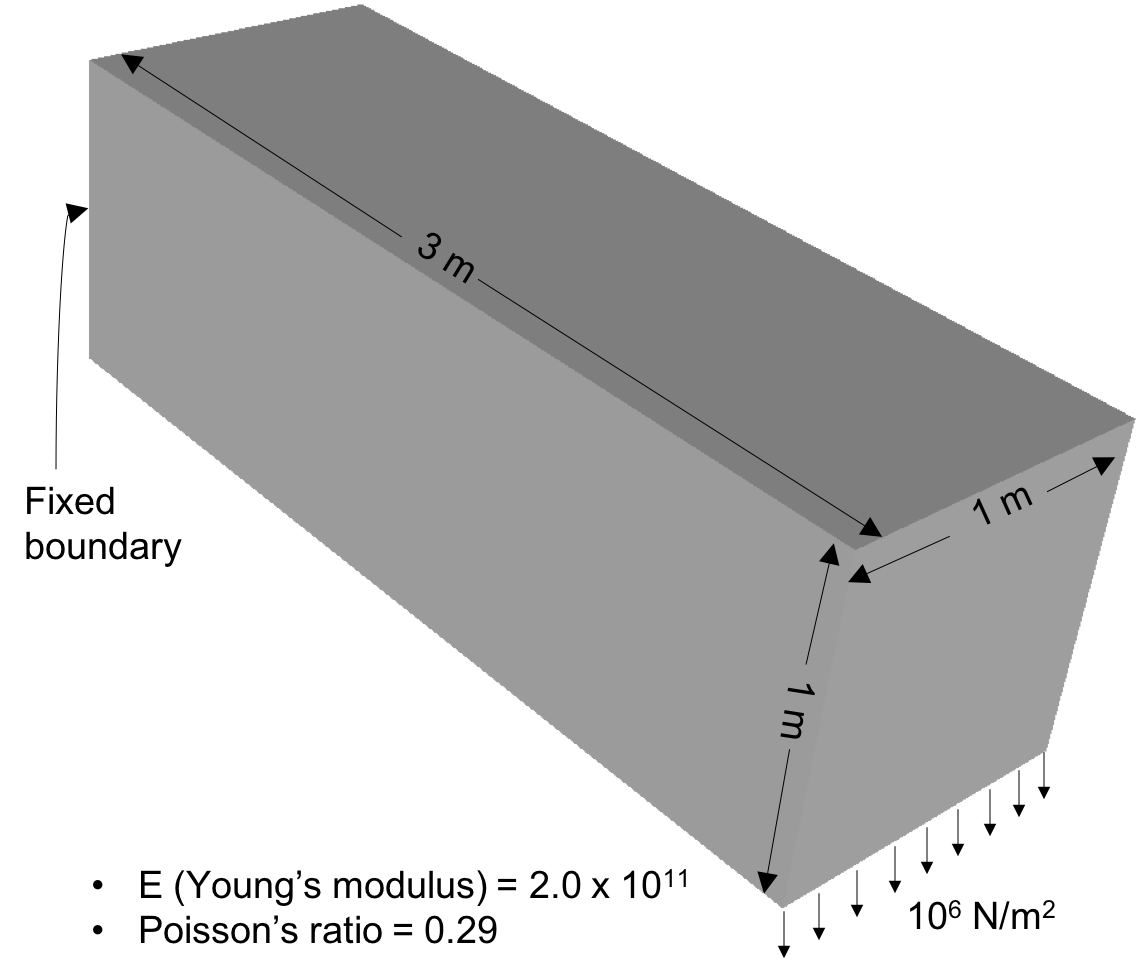}
      \quad\quad\quad
    \includegraphics[scale=0.18]{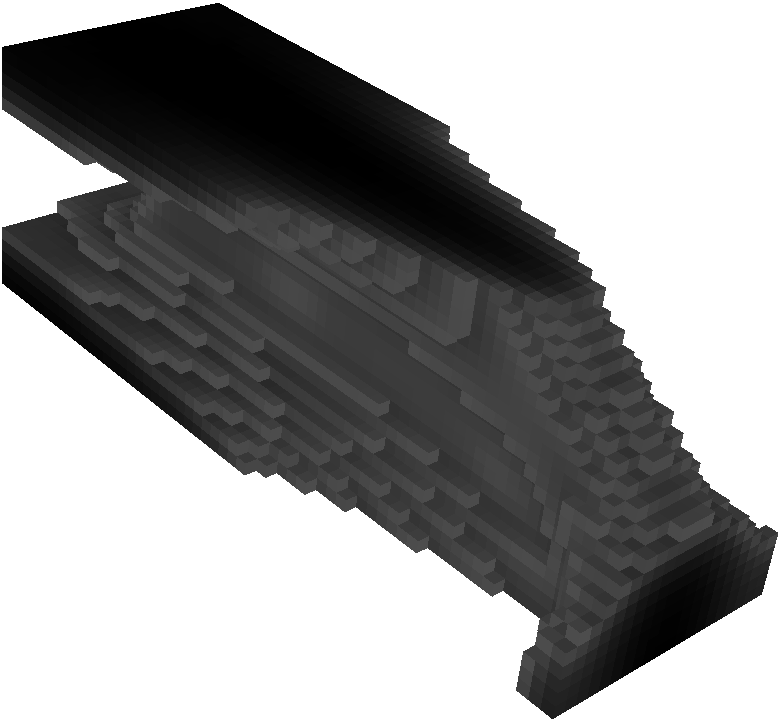}
    \end{center}
    \caption{Left: design domain, Right: an optimal design for the cantilever
    beam}
    \label{fig:design3Dbeam}
  \end{figure}

  Cantilever beam design problem is the same as the numerical examples
  considered in \cite{wang2007large}. The design domains and optimal designs for
  two problems are shown in Figure~\ref{fig:design3Dbeam}. Compliance
  minimization with total mass constraint is considered.
  Three different mesh resolutions are considered: $3\ 000$, $24\ 000$, and $192\
  000$ design variables.  
  The following material properties are used: Young's modulus of
  $2.0\times10^11$ $N/m^2$ and Poisson's ratio of $0.29$.  The upper bound for
  the mass constraint is $0.5$ and the SIMP parameter is $\penaltySIMP = 3$ in
  Eq.~\eqref{eq:SIMP}.  IPOPT is used as an optimization solver with convergence
  threshold of $10^{-6}$.  We use the Helmholtz filter with $\helmholtzradius =
  0.1\ m$. Finally, the following ROM-based topology optimization
  parameters are used: $\scaleROM = \scalePCG = 10^{-2}$, $\scaleCUT =
  10^{-3}$, $\absthresholdPCG = 10^{-4}$, $\qrtol = \svdtol = 10^{-9}$, and
  $\maxrank = 10$. All the simulations for the cantilever beam design problem
  use $36$ processors in Quartz. 
  A structured mesh with uniform hexahedral
  first-order finite elements is used for the discretization.

  \begin{table}[h]\begin{center}
  \begin{tabular}{l*{3}{r}}
      & \bf Default & \bf ROM-based top.opt.  & \bf ROM-based top.opt.  \\ 
      &             & \bf incremental QR   & \bf incremental SVD  \\
    \hline
    Optimal compliance          & 13.59   & 13.59   & 13.63  \\
    Total wall clock time of linear solve (sec.) & \bad{41.9} & \good{26.2} &
    \good{10.2} \\
    Total iters. of linear solve & \bad{7,700} & \good{1,528} & \good{500} \\
    Avg. iters. of linear solve & \bad{23.5} & \good{5.2} & \good{4.3} \\
    IPOPT iter.                 & 297 & 294 & 113 \\
    KKT norm                    & 6.68e-7 & 4.86e-7 & 6.56e-7 \\
    \hline
    Speed-up of total linear solve        & 1.0 & \good{1.6}  & \good{4.1}  \\
    Avg. iter. reduction of linear solve  & 1.0 & \good{2.4} & \good{5.5} \\
    Total iter. reduction of linear solve & 1.0 & \good{5.0} & \good{15.4}
  \end{tabular}\end{center}
    \caption{performance comparison for cantilever beam 3D design problem with
    $3,000$ design variables.}
  \label{ta:3Dcantilever0.1}
  \end{table}

  \begin{table}[h]\begin{center}
  \begin{tabular}{l*{3}{r}}
      & \bf Default & \bf ROM-based top.opt.  & \bf ROM-based top.opt.  \\ 
      &             & \bf incremental QR   & \bf incremental SVD  \\
    \hline
    Optimal compliance          &  13.2  & 13.2   & 13.2  \\
    Total wall clock time of linear solve (sec.) & \bad{241.2} & \good{108.4} &
    \good{69.8} \\
    Total iters. of linear solve & \bad{6,789} & \good{1,083} & \good{799} \\
    Avg. iters. of linear solve & \bad{24.2} & \good{6.4} & \good{7.3} \\
    IPOPT iter.                 & 180 & 168 & 108 \\
    KKT norm                    & 7.74e-7 & 6.73e-7 & 8.11e-7 \\
    \hline
    Speed-up of total linear solve        & 1.0 & \good{2.2} & \good{3.5}  \\
    Avg. iter. reduction of linear solve  & 1.0 & \good{3.8} & \good{3.3} \\
    Total iter. reduction of linear solve & 1.0 & \good{6.3} & \good{8.5}
  \end{tabular}\end{center}
    \caption{performance comparison for cantilever beam 3D design problem with
    $24,000$ design variables.}
  \label{ta:3Dcantilever0.05}
  \end{table}

  \begin{table}[h]\begin{center}
  \begin{tabular}{l*{3}{r}}
      & \bf Default & \bf ROM-based top.opt.  & \bf ROM-based top.opt.  \\ 
      &             & \bf incremental QR   & \bf incremental SVD  \\
    \hline
    Optimal compliance          & 13.6 & 13.6  & 13.6  \\
    Total wall clock time of linear solve (sec.) & \bad{761.8} & \good{267.5} &
    \good{238.3} \\
    Total iters. of linear solve & \bad{3,229} & \good{942} & \good{989} \\
    Avg. iters. of linear solve & \bad{24.1} & \good{6.5} & \good{7.9} \\
    IPOPT iter.                 & 132 & 142 & 123 \\
    KKT norm                    & 9.5-e-7 & 8.07e-7 & 9.86e-7 \\
    \hline
    Speed-up of total linear solve        & 1.0 & \good{2.8} & \good{3.2}  \\
    Avg. iter. reduction of linear solve  & 1.0 & \good{3.7} & \good{3.1} \\
    Total iter. reduction of linear solve & 1.0 & \good{3.4} & \good{3.3}
  \end{tabular}\end{center}
    \caption{performance comparison for cantilever beam 3D design problems with
    $192,000$ design variables.}
  \label{ta:3Dcantilever0.025}
  \end{table}

  Tables~\ref{ta:3Dcantilever0.1}, \ref{ta:3Dcantilever0.05}, and
  \ref{ta:3Dcantilever0.025} show the performance comparison among the default
  and ROM-based topology optimization methods with incremental QR and SVD for
  the 3D cantilever beam design problems with different number of design
  variables. They have the comparable optimal compliance values although
  ROM-based approach with incremental SVD gives slightly worse optimal value
  with $3,000$ design variables in Table~\ref{ta:3Dcantilever0.1}. It
  converges to a different local minimum. The ROM-based approach with
  incremental SVD converges within $113$ IPOPT iterations, while the other
  approaches take around $300$ iterations. Therefore, it results in a tremendous
  reduction in total number of iterations in linear solve, i.e., $15.4$ as well
  as a considerable speed-up for the solving time, i.e., $4.1$. 
  
  Other cases with more design
  variables, i.e., Tables~\ref{ta:3Dcantilever0.05} and
  \ref{ta:3Dcantilever0.025}, all the methods converge to the same design with
  the same optimal compliance value.  In terms of the total wall clock time of
  linear solves, the ROM-based approach achieves a considerable speed-up of
  $2.2$ to $3.5$.  The average number of iterations in linear solves are also
  reduced considerably by a factor of up to $3.1$ to $3.8$. For these particular
  problems, the ROM-based approaches converge to optimal solutions with less
  number of IPOPT iterations than the default method. Therefore, the reduction
  factor of the total number linear solve iterations is quite big, i.e., from
  $3.3$ to $8.5$.

  \subsubsection{Wind turbine blade design}\label{sec:bladedesign}
\begin{figure}[th]
    \begin{center}
    \includegraphics[scale=0.35]{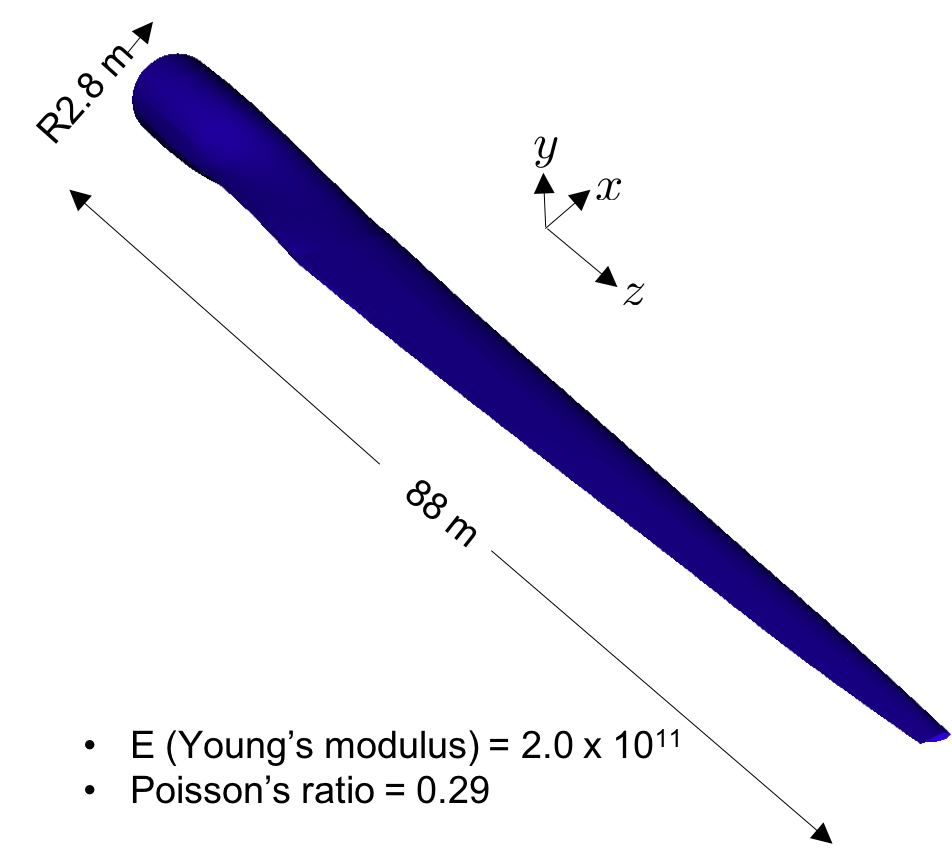}
      \quad\quad\quad
    \includegraphics[scale=0.20]{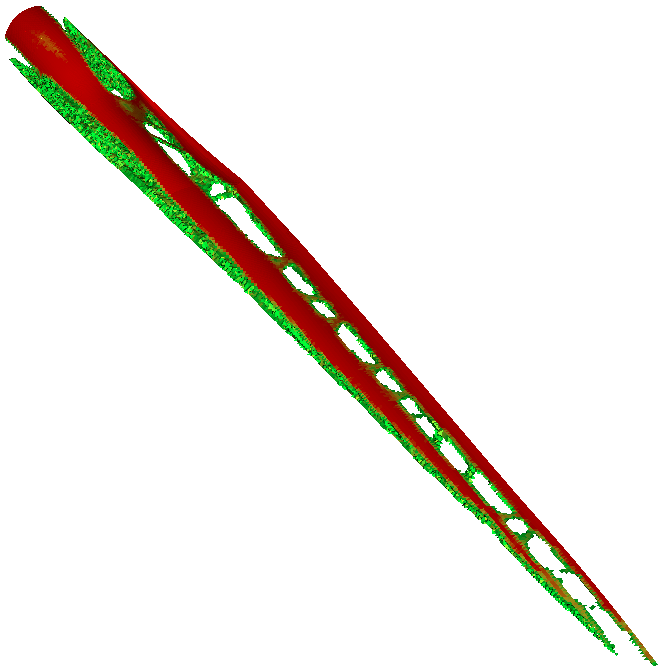}
    \end{center}
    \caption{Left: design domain, Right: an optimal design for the
    wind turbine blade}
    \label{fig:design3Dblade}
  \end{figure}

  To obtain an optimal blade design problem, we consider minimizing compliance
  with a total mass constraint.  The design domain is described in
  Figure~\ref{fig:design3Dblade}.  A fixed boundary condition is applied to the
  thicker end of the blade.  The blade problem has $414,979$ design variables.
  In order to mimic wind conditions, two load cases are considered: $10$ $N$ in
  $x$-direction and $10$ $N$ in $y$-direction.  Two different compliance values
  are computed for two different loads, then the average compliance is
  minimized.  The following material properties are used: Young's modulus of
  $2.0\times10^11$ $N/m^2$ and Poisson's ratio of $0.29$.  The upper bound for
  the mass constraint is $0.25$. The SIMP parameter $\penaltySIMP = 3$ is used
  in Eq.~\eqref{eq:SIMP}.  IPOPT is used as an optimization solver with
  convergence threshold of $10^{-6}$.  The mass matrix filter described in
  Section~\ref{sec:topopt} is used to avoid checkerboard problem.  Finally, the
  following ROM-based topology optimization parameters are used: $\scaleROM =
  \scalePCG = 10^{-3}$, $\scaleCUT = 10^{-3}$, $\absthresholdPCG = 10^{-4}$,
  $\qrtol = \svdtol = 10^{-9}$, and $\maxrank = 10$. All the simulations for the
  blade design problem use $72$ processors from Quartz.

  \begin{table}[h]\begin{center}
  \begin{tabular}{l*{3}{r}}
      & \bf Default & \bf ROM-based top.opt.  & \bf ROM-based top.opt.  \\ 
      &             & \bf incremental QR   & \bf incremental SVD  \\
    \hline
    Optimal compliance        & 2.61 & 2.61   & 2.61  \\
    Total wall clock time of linear solve (hour) & \bad{1.70} & \good{0.55} &
    \good{0.48} \\
    Total iters. of linear solve & \bad{372,196} & \good{56,542} &
    \good{51,062} \\
    Avg. iters. of linear solve & \bad{265.5} & \good{25.8} & \good{27.3} \\
    IPOPT iter. & 664 & 1092 & 934 \\
    KKT norm & 9.82e-7 & 8.71e-7 & 7.47e-7 \\
    \hline
    Speed-up of total linear solve        & 1.0 & \good{3.1}  & \good{3.5}  \\
    Avg. iter. reduction of linear solve  & 1.0 & \good{10.3} & \good{9.5} \\
    Total iter. reduction of linear solve & 1.0 & \good{6.6} & \good{9.7}
  \end{tabular}\end{center}
    \caption{performance comparison for wind turbine blade 3D design problem.}
  \label{ta:optimalblade-results}
  \end{table}

  One difference between the blade problem and the 3D cantilever beam problme in
  the previous section is mesh.  An unstructured mesh with tetrahedral
  first-order finite elements is used for the discretization of the blade.
  Because of the unstructured mesh, the default method with AMG preconditioner
  from HYPRE is not enough to bring down the number of PCG iterations as shown
  in Table~\ref{ta:optimalblade-results}, i.e., $265.5$ average number of
  iterations per a linear solve is required for the default method. On the other
  hand, the ROM-based approaches take only $25.8$ and $27.3$ average number of
  iterations for incremental QR and SVD, respectively.  This gives the reduction
  factor of around $10$.  Our method also achieves considerable speed-ups in
  terms of linear system solving times, i.e., larger than $3$, reducing total
  wall-clock time from $1.7$ hours to $0.55$ or $0.48$ hours, even though our
  ROM-based approaches take more IPOPT iterations than the default method.  Note
  that the three different methods produce the identical optimal compliance
  value, i.e., $2.61$. These optimal designs satisfy the KKT optimality
  conditions, implying that the quality of the design is not degraded by the
  approximation introduced by the ROMs.

  \subsection{Stress-constrained problem}\label{sec:stress}
  Our method can be applicable not only to the compliance minimization problems,
  but also to stress-constrained design problems. We demonstrate it in this
  section by considering a classical stress constrained topology optimization
  problem, i.e., L-bracket problem. 

  \subsubsection{L-bracket problem}\label{sec:lbracket}
\begin{figure}[th]
    \begin{center}
    \includegraphics[scale=0.35]{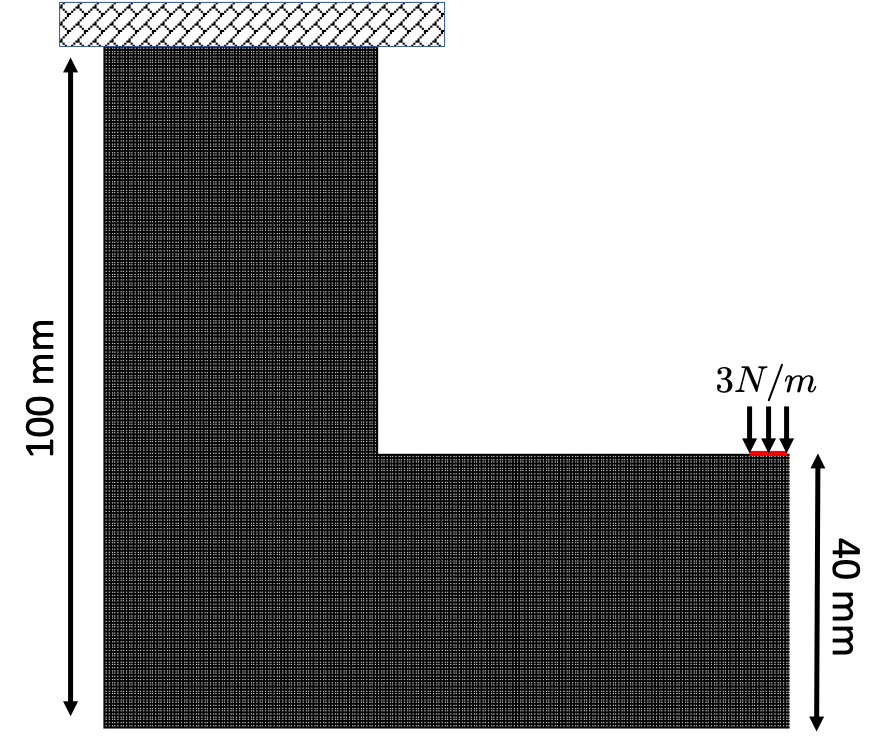}\quad\quad
    \includegraphics[scale=0.16]{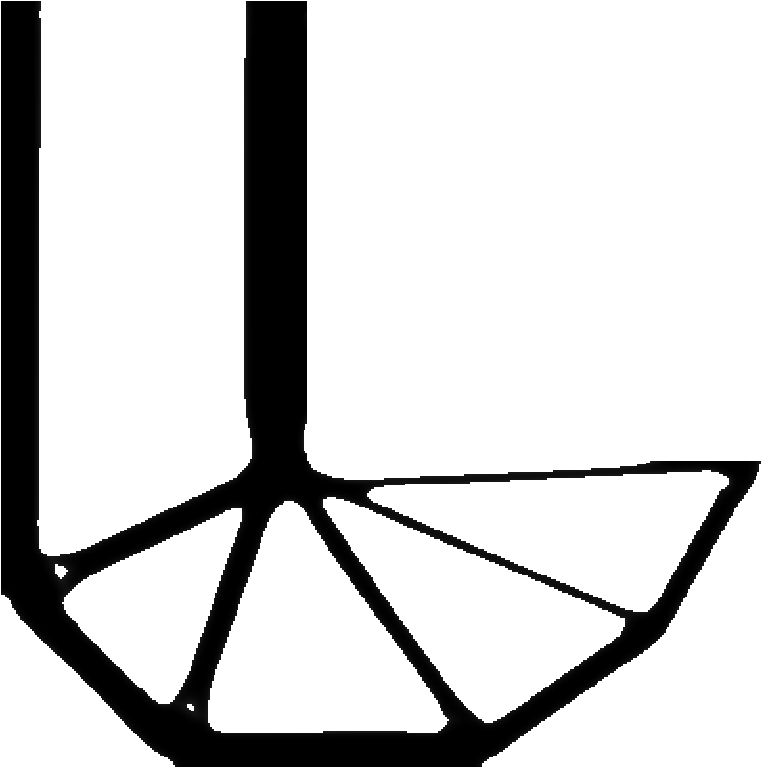}
    \end{center}
    \caption{L-bracket stress-constrained problem, Left: design domain, boundary
    conditions, and external load, Right: an optimal design}
    \label{fig:lbracket_design}
  \end{figure}

  A classical L-bracket stress-constrained problem is considered.  The design
  domain, boundary conditions, and external loading are described in
  Figure~\ref{fig:lbracket_design}.  The total mass is minimized with a stress
  constraint.  The von Mises stress criterion is used as a stress quantity and
  p-norm with $p=8$ is used to approximate the maximum stress value as in
  \cite{le2010stress}.  The following material properties are used: Young's
  modulus of $10^6$ $N/m^2$ and Poisson's ratio of $0.3$.  The number of design
  variables is $102,400$.  The upper bound for the stress constraint is $20
  N/m^2$ and the SIMP parameter is $\penaltySIMP = 3$ in Eq.~\eqref{eq:SIMP}.
  Stress quantity penalization parameter $\penaltyStress = 0.5$ in
  Eq.~\eqref{eq:relaxedStress} is used. The Helmholtz filter with
  $\helmholtzradius = 0.0005\ m$ is used.  IPOPT is used as an optimization solver with
  convergence threshold of $10^{-6}$.  Finally, we use the following ROM-based
  topology optimization parameters: $\scaleROM = \scalePCG = 10^{-3}$,
  $\scaleCUT = 10^{-3}$, $\absthresholdPCG = 10^{-4}$, $\qrtol = \svdtol =
  10^{-9}$, and $\maxrank = 10$. All the simulations for the stress-constrained
  problems use $144$ processors in Quartz.

  \begin{table}[h]\begin{center}
  \begin{tabular}{l*{3}{r}}
      & \bf Default & \bf ROM-based top.opt.  & \bf ROM-based top.opt.  \\ 
      &             & \bf incremental QR   & \bf incremental SVD  \\
    \hline
    Optimal mass             & 0.38 & 0.38 & 0.38  \\
    Total wall clock time of linear solve (hour) & \bad{1.5} & \good{0.7} &
    \good{0.71} \\
    Total iters. of linear solve & \bad{721,223} & \good{290,245} &
    \good{325,418} \\
    Avg. iters. of linear solve & \bad{61.0} & \good{39.5} & \good{46.5} \\
    IPOPT iter. & 4,788 & 3,670 & 3,494 \\
    KKT norm & 7.78e-7 & 6.00e-7 & 9.83e-7 \\
    \hline
    Speed-up of total linear solve        & 1.0 & \good{2.1} & \good{2.1}  \\
    Avg. iter. reduction of linear solve  & 1.0 & \good{1.5} & \good{1.3} \\
    Total iter. reduction of linear solve & 1.0 & \good{2.5} & \good{2.2}
  \end{tabular}\end{center}
    \caption{performance comparison for L-bracket stress-constrained problem.}
  \label{ta:lbracket-results}
  \end{table}

  A structured mesh with uniform quadrilateral first-order finite elements is
  used for the discretization.  Therefore, the default method with AMG
  preconditioner is able to reduce the number of linear solve iteration
  sufficiently. However, we still see a further reduction and speed-up by
  applying our ROM-based approach.  For example, the default method requires
  $61.0$ iterations in average per a linear solve, while the ROM-based
  approaches with incremental QR and SVD require $39.5$ and $46.5$ iterations,
  respectively. This gives reduction of $1.5$ and $1.3$, respectively. The wall
  clock time for the default method is $1.5$ hours, while the ROM-based
  approaches finish within $0.71$ hours, resulting in a speed-up of $2.1$.  Note
  that all the three methods produce the same optimal mass, which is $0.38$.
  Also note that all the three methods converge to a point that satisfies the
  necessary optimality conditions, i.e., the KKT norms are less than $1.0e-6$.

\section{Conclusion}\label{sec:conclusion}
  A ROM-based design optimization acceleration method is introduced.  The
  overall design optimization is accelerated by accelerating linear system
  solves as demonstrated in Section~\ref{sec:numeric}. The ROM-based approach
  shows a considerable speed-up especially when the unstructured mesh is used,
  in which the default method with a AMG preconditioner requires many Krylov
  subspace iterations.  The method is not tailored for the compliance
  minimization problems.  It is applicable to a stress-constrained optimization
  problem, which is also demonstrated in numerical experiments. Furthermore, our
  method is general enough to be applicable to other PDE-constrained
  optimization problems, such as shape optimization and inverse problems.
  Finally, the method does not suffer from the approximation introduced by the
  ROM because the accuracy of ROM is carefully monitored and treated throughout
  the optimization process, resulting in an optimal design that satisfies the
  KKT optimality condition. 

  Future research is required to further understand the precise conditions for
  the inexactness. We have only provided the heuristic explanations on why our
  ROM-based approach works well and determines the parameter values of our
  method heuristically. Thorough theoretical study on the convergence rate of
  the interior-point method affected by the inexactness introduced by ROMs is
  necessary because the majority of literatures on this topic considers the
  inexactness coming from the optimization linear solves, not from the PDE
  solves.  Finally, constructing reduced order operator can be computationally
  expensive. In reduced order model research community, a hyper-reduction is
  used to reduce the cost of constructing reduced order operators.  This will be
  investigated in future to further accelerate the optimization process.

\section*{Acknowledgments}
This work was performed at Lawrence Livermore National Laboratory and was
supported by the LDRD program (17-SI-005). 
Lawrence Livermore National Laboratory is operated by Lawrence
Livermore National Security, LLC, for the U.S. Department of Energy,
National Nuclear Security Administration under Contract DE-AC52-07NA27344
and LLNL-JRNL-791183.

  \bibliographystyle{plain}
  \bibliography{references}

\end{document}